\documentclass[smallextended]{svjour3}       
\smartqed  
\usepackage{graphicx} 
\usepackage[font=small,labelfont=bf]{caption} 
\usepackage[subrefformat=parens]{subcaption}

\usepackage[cp1251]{inputenc}

\usepackage[dvipsnames]{xcolor}
\definecolor{br}{RGB}{201,115,6}
\definecolor{gr}{RGB}{11,156,19}
\definecolor{rd}{RGB}{225,18,124}
\definecolor{brl}{RGB}{246,181,0}
\definecolor{grl}{RGB}{7,157,2}
\definecolor{rdl}{RGB}{180,11,77}

\usepackage{color}
\usepackage[colorlinks]{hyperref}
\usepackage{mathptmx}      
\usepackage{amsmath,amssymb,amsbsy,euscript,amsfonts,cite}
\begin{document}

\title{Approximation of Functions of Several
Variables by Multidimensional \emph{A}- and \emph{J}-fractions
with Independent Variables}

\titlerunning{Approximation of Functions of Several
Variables by Branched Continued Fractions}        

\author{Roman~Dmytryshyn \and Serhii~Sharyn}

\institute{Roman~Dmytryshyn \and Serhii~Sharyn\at
    Department of Mathematical and Functional Analysis,
    Vasyl Stefanyk Precarpathian National University,
    str. Shevchenka~57, Ivano-Frankivsk~76018, Ukraine\\
    \email{dmytryshynr@hotmail.com, serhii.sharyn@pnu.edu.ua}
    }

\date{\today}

\maketitle

\begin{abstract}
The paper deals with the problem of approximating the functions of
several variables by branched continued fractions, in particular,
multidimensional \emph{A}- and \emph{J}-fractions with independent
variables. A generalization of Gragg's algorithm is constructed
that enables us to compute, by the coefficients of the given
formal multiple power series, the coefficients of the
corresponding multidimensional \emph{A}-fraction with independent
variables. This algorithm can also be used to construct the
multidimensional \emph{J}-fraction with independent variables
corresponding to a given formal multiple Laurent series. Some
numerical experiments of approximating the functions of several
variables by these branched continued fractions are given.
\end{abstract}

\section{Introduction}
The problem of representing the functions of several variables,
which arises, in particular, when solving different functional
equations, contributes to the development and implementation of
effective methods and algorithms that are implemented until the
construction of special software. Currently, various tools are
used to represent and/or approximate such functions. Possibly one
of the most effective is the multidimensional generalization of
continued fractions -- branched continued fractions \cite{bod86}.
The construction of the rational approximations of a function of
several variables is based on a correspondence between the
approximants of the branched continued fraction and the formal
multiple power series, which represents this function (see,
\cite{bodm19,dm20prsesa}). Furthermore, the problem of
constructing the corresponding branched continued fractions
contributed to the emergence of their various structures (see, for
example, \cite{bod98,bod03,cuve88}).

In \cite{bod76}, D. Bodnar introduced so-called ``branched
continued fractions with independent variables'', which by their
structure are a multidimensional analogue of multiple power
series. The correspondence properties of these branched continued
fractions with polynomial elements are closely connected to the
degree and form of these polynomials. A number of their types are
essential in expanding the analytic functions of several variables
in terms of branched continued fractions with independent
variables. Using the idea of constructing an algorithm in
\cite{bod98}, in \cite{bar14} O. Baran constructed
multidimensional \emph{C}- and \emph{A}-fractions with independent
variables corresponding to a given formal multiple power series.

In \cite{dm12jat}, the author proposed an algorithm based on the
classical algorithm that enables us to compute, by the
coefficients of the given formal double power series, the
coefficients of the corresponding two-dimensional
\emph{g}-fractions with independent variables. This idea was used
to construct the corresponding multidimensional \emph{A}- and
\emph{J}-fractions with independent variables in \cite{bodm19},
multidimensional regular \emph{C}-fractions with independent
variables in \cite{dm20prsesa}, and multidimensional
\emph{S}-fractions with independent variables in \cite{dmsh21}. On
the one hand, a few numerical experiments show the proposed
algorithm's efficiency and, on the other hand, the power and
feasibility of the method to numerically approximate analytic
functions of several variables from their formal multiple power
series (see more examples in \cite{dm12jms}).

In \cite{rut54a,rut54b,rut54c} H. Rutishauser proposed a
quotient-difference algorithm (\emph{qd}-algo\-rithm) that
provides a convenient numerical procedure for computing the
coefficients of continued fractions corresponding to the given
formal power series. Various generalizations of the
\emph{qd}-algorithm to the problems of multivariate rational
interpolation and multivariate Pad\'e approximants were developed
and applied by G. Claessens \cite{cl81} and A. Cuyt
\cite{cuyt88,cuyt83}. In \cite{dm15jms}, the author proposed a
two-dimensional generalization of \emph{qd}-algorithm for
computing the coefficients of the corresponding two-dimensional
regular \emph{C}-fractions with independent variables.

In this paper, we construct a generalization of Gragg's algorithm
based on the algorithm given in \cite{gr74}, which allows us to
compute, by the coefficients of the given formal multiple power
series (formal multiple Laurent series), the coefficients of the
corresponding multidimensional \emph{A}-fraction with independent
variables (multidimensional \emph{J}-fraction with independent
variables). As an application, we provide some numerical
experiments.

\section{Preliminaries}\label{sec:2}
Let $N$ be a fixed natural number, $\mathbb{Z}_{\ge0}$ be the set
of non-negative integers, $\mathbb{C}$ be the set of complex
numbers,
$\mathbb{Z}^N_{\ge0}=\mathbb{Z}_{\ge0}\times\mathbb{Z}_{\ge0}\times\ldots\times\mathbb{Z}_{\ge0}$
be the Cartesian product of $N$ copies of the $\mathbb{Z}_{\ge0},$
$\mathbb{C}^N=\mathbb{C}\times\mathbb{C}\times\ldots\times\mathbb{C}$
be the Cartesian product of $N$ copies of the $\mathbb{C},$
${\bf{k}}=(k_1,k_2,\ldots,k_N)$ be an element of
$\mathbb{Z}^N_{\ge0},$ ${\bf{z}}=(z_1,z_2,\ldots,z_N)$ be an
element of $\mathbb{C}^N.$ For ${\bf{k}}\in\mathbb{Z}^N_{\ge0}$
and ${\bf{z}}\in\mathbb{C}^N$ put
\begin{align*}
{\bf{k}}!=k_1!k_2!\ldots
k_N!,\quad\vert{\bf{k}}\vert=k_1+k_2+\ldots+k_N,\quad{\bf{z}}^{{\bf{k}}}=z_1^{k_1}z_2^{k_2}\ldots
z_N^{k_N}.
\end{align*}

A series of the form
\begin{align*}
L({\bf{z}})=\sum_{\vert{\bf{k}}\vert\ge{0}}c_{{\bf{k}}}{\bf{z}}^{{\bf{k}}},
\end{align*}
where $c_{{\bf{k}}}\in\mathbb{C}$ for $\vert{\bf{k}}\vert\ge0,$ is
called a formal multiple power series at ${\bf{z}}={\bf0}.$ A set
of formal multiple power series at ${\bf{z}}={\bf0}$ is denoted by
$\mathbb{L}.$

Let~$R({\bf{z}})$ be a function holomorphic in a neighbourhood of
the origin $({\bf z}={\bf0}).$ Let the mapping $\Lambda:R({\bf
z})\to\Lambda(R)$ associates with $R({\bf z})$ its Taylor
expansion in a neighbourhood of the~origin. A sequence $\{R_n({\bf
z})\}$ of functions holomorphic at the origin is said to
correspond at ${\bf z}={\bf0}$ to a formal multiple power series
$L({\bf z})$ if
\begin{align*}
\lim_{n\to\infty}\lambda\big(L-\Lambda(R_n)\big)=\infty,
\end{align*}
where $\lambda$ is the function defined as follows:
$\lambda:\mathbb{L}\to\mathbb{Z}_{\ge0}\cup\{\infty\}$; if $L({\bf
z})\equiv0$ then $\lambda(L)=\infty;$ if $L({\bf z})\not\equiv0$
then $\lambda(L)=m,$ where $m$ is the smallest degree of
homogeneous terms for which $c_{{\bf{k}}}\ne0,$ that is
$m=\vert{\bf{k}}\vert.$

If $\{R_n({\bf z})\}$ corresponds at ${\bf z}={\bf0}$ to a formal
multiple power series $L({\bf z}),$ then the order of
correspondence of $R_n({\bf z})$ is defined to be
\begin{align*}
\nu_n=\lambda(L-\Lambda(R_n)).
\end{align*}
By the definition of $\lambda,$ the series $L({\bf{z}})$ and
$\Lambda(R_n)$ agree for all homogeneous terms up to and including
degree $(\nu_n-1).$

Let $i(0)=0$ and $\mathcal{I}_0=\{0\}.$ Let us introduce the
following sets of multiindices for $k\ge1$
\begin{align*}
\mathcal{I}_k=\{i(k):\;i(k)=(i_1,i_2,\ldots,i_k),\;1\le{i}_p\le{i}_{p-1},\;1\le{p}\le{k},\;i_0=N\}.
\end{align*}
Let
$\langle\{a_{i(k)}\}_{i(k)\in\mathcal{I}_k,\;k\ge1},\{b_{i(k)}\}_{i(k)\in\mathcal{I}_k,\;k\ge0}\rangle$
denote the ordered pair of sequences of complex numbers with
$a_{i(k)}\ne0$ for all $i(k)\in\mathcal{I}_k,$ $k\ge1,$ and if for
$k\ge1$ there exist a multiindex $i(k)\in\mathcal{I}_k$ such that
$b_{i(k)}=0,$ than $b_{i(k-1),j}\ne0$ for $1\le j\le i_{k-1}$ and
$j\ne i_k.$ Let the sequence $\{f_k\}$ is defined as follows:
\begin{align*}
f_0&=b_0,\\
f_1&=b_0+\sum_{i_1=1}^N\frac{a_{i(1)}}{b_{i(1)}},\\
f_2&=b_0+\sum_{i_1=1}^N\frac{a_{i(1)}}{b_{i(1)}}{\atop+}\sum_{i_2=1}^{i_1}\frac{a_{i(2)}}{b_{i(2)}},
\\&\ldots,\\
f_k&=b_0+\sum_{i_1=1}^N\frac{a_{i(1)}}{b_{i(1)}}{\atop+}\sum_{i_2=1}^{i_1}\frac{a_{i(2)}}{b_{i(2)}}{\atop+}\ldots{\atop+}
\sum_{i_k=1}^{i_{k-1}}\frac{a_{i(k)}}{b_{i(k)}},\\&\ldots.
\end{align*}

The ordered pair
\[\langle\langle\{a_{i(k)}\}_{i(k)\in\mathcal{I}_k,\;k\ge1},\{b_{i(k)}\}_{i(k)\in\mathcal{I}_k,\;k\ge0}\rangle,
\{f_k\}_{k\ge0}\rangle\] is the branched continued fraction with
independent variables denoted by the symbol
\begin{equation}\label{eq:1}
b_0+\sum_{i_1=1}^N\frac{a_{i(1)}}{b_{i(1)}}{\atop+}\sum_{i_2=1}^{i_1}\frac{a_{i(2)}}{b_{i(2)}}{\atop+}
\sum_{i_3=1}^{i_2}\frac{a_{i(3)}}{b_{i(3)}}{\atop+}\ldots
\end{equation}
(note that here we used the more convenient notation of
(\ref{eq:1}), proposed by J.F.W.~Herschel (see,
\cite[p.~148]{hers20}) following the example of H.H.~B\"urmann).
The numbers $a_{i(k)}$ and $b_{i(k)}$ are called the elements of
the branched continued fraction with independent variables, the
relation $a_{i(k)}/b_{i(k)}$ is called the $k$th partial quotient
and the value $f_k$ is called the $k$th approximant.

Let $(i_1,i_2,\ldots,i_k,\ldots)$ be a fixed infinite multiindex
such that $1\le i_k\le i_{k-1}$ for $k\ge1,$ where $i_0=N.$ The
continued fraction
\begin{align*}
\frac{a_{i_1}}{b_{i_1}}{\atop+}\frac{a_{i_1,i_2}}{b_{i_1,i_2}}{\atop+}\frac{a_{i_1,i_2,i_3}}{b_{i_1,i_2,i_3}}{\atop+}\ldots
\end{align*}
is called the $(i_1,i_2,\ldots,i_k,\ldots)$-branch of the branched
continued fraction with independent variables (\ref{eq:1}).

Next, let $e_0=(0,0,\ldots,0),$
$e_k=(\delta_{k,1},\delta_{k,2},\ldots,\delta_{k,N})$ be a
multiindex, where $1\le k\le N,$ $\delta_{i,j}$ is a Kronecker
symbol. Let us introduce the following sets of multiindices for
$k\ge1$
\begin{align*}
\mathcal{E}_k=\{e_{i(k)}:\;e_{i(k)}=e_{i_1,i_2,\ldots,i_k}=e_{i_1}+e_{i_2}+\ldots+e_{i_k},\;i(k)\in\mathcal{I}_k\}
\end{align*}
and the mapping $\varphi:$ $\mathcal{I}_k\to\mathcal{E}_k,$ such
that $\varphi(i(k))=e_{i(k)}$ for all $i(k)\in\mathcal{I}_k,$
$k\ge1.$ It can be shown that the mapping $\varphi$ is bijective.

A branched continued fraction with independent variables of the
form
\begin{equation}\label{eq:2}
\sum_{i_1=1}^N\frac{p_{e_{i(1)}}z_{i_1}}{1+q_{e_{i(1)}}z_{i_1}}{\atop+}
\sum_{i_2=1}^{i_1}\frac{(-1)^{\delta_{i_1,i_2}}p_{e_{i(2)}}z_{i_1}z_{i_2}}{1+q_{e_{i(2)}}z_{i_2}}{\atop+}
\sum_{i_3=1}^{i_2}\frac{(-1)^{\delta_{i_2,i_3}}p_{e_{i(3)}}z_{i_2}z_{i_3}}{1+q_{e_{i(3)}}z_{i_3}}{\atop+}\cdots,
\end{equation}
where the $p_{e_{i(k)}}\in\mathbb{C}\setminus\{0\},$
$q_{e_{i(k)}}\in\mathbb{C}$ for $e_{i(k)}\in\mathcal{E}_k,$
$k\ge1,$ is called a multidimensional \emph{A}-fraction with
independent variables. The $n$th approximant $f_n({\bf{z}})$ of
(\ref{eq:2}) is expressed by
\begin{align*}
f_n({\bf{z}})&=\sum_{i_1=1}^N\frac{p_{e_{i(1)}}z_{i_1}}{1+q_{e_{i(1)}}z_{i_1}}\\&\quad{\atop+}
\sum_{i_2=1}^{i_1}\frac{(-1)^{\delta_{i_1,i_2}}p_{e_{i(2)}}z_{i_1}z_{i_2}}{1+q_{e_{i(2)}}z_{i_2}}{\atop+}\ldots{\atop+}
\sum_{i_n=1}^{i_{n-1}}\frac{(-1)^{\delta_{i_{n-1},i_n}}p_{e_{i(n)}}z_{i_{n-1}}z_{i_n}}{1+q_{e_{i(n)}}z_{i_n}},\quad
n\ge1.
\end{align*}

A multidimensional \emph{A}-fraction with independent variables
(\ref{eq:2}) is said to correspond at ${\bf z}={\bf0}$ to a formal
multiple power series $L({\bf z})$ if its sequence of approximants
$\{f_n({\bf{z}})\}$ corresponds to $L({\bf z})$ at ${\bf
z}={\bf0}.$

\section{Generalization of Gragg's algorithm}\label{sec:3}
Let $N\ge2.$ We will construct an algorithm for the expansion of
the given formal multiple power series
\begin{equation}\label{eq:3}
L({\bf{z}})=\sum_{\vert{\bf{k}}\vert\ge1}c_{{\bf{k}}}{\bf{z}}^{{\bf{k}}},
\end{equation}
where $c_{{\bf{k}}}\in\mathbb{C}$ for $\vert{\bf{k}}\vert\ge1,$
into the corresponding multidimensional \emph{A}-fraction with
independent variables (\ref{eq:2}) in much the same way as the
algorithm in \cite{bodm19}. The process of constructing will be
shown step by step.

Step 1.1:  Let $c_{e_{i_1}}\ne0$ for $2\le i_1\le N.$ Then we can
rewrite $L(\bf{z})$ in the form
\begin{align*}
L({\bf{z}})=P_{e_0}(z_1)+\sum_{{i_1}=2}^Nc_{e_{i_1}}z_{i_1}R_{e_{i_1}}({\bf{z}}),
\end{align*}
where
\begin{align*}
P_{e_0}(z_1)=\sum_{n=1}^\infty{c}_{ne_1}z_1^n,\quad
R_{e_{i_1}}({\bf{z}})=\sum_{\substack{{\bf{k}}\ge{\bf0}\\k_j=0,\;
i_1+1\le{j}\le{N}}}\frac{c_{{\bf{k}}+e_{i_1}}}{c_{e_{i_1}}}{\bf{z}}^{{\bf{k}}}.
\end{align*}

Step 1.2:  Let $H_{e_1}(n)\ne0$ for $n\ge1,$ where
\begin{align*}
H_{e_1}(n)=\left\vert\begin{array}{cccc}
c_{e_1}&c_{2e_1}&\ldots&c_{ne_1}\\
c_{2e_1}&c_{3e_1}&\ldots&c_{(n+1)e_1}\\\hdotsfor{4}\\
c_{ne_1}&c_{(n+1)e_1}&\ldots&c_{(2n-1)e_1}\end{array}\right\vert
\end{align*}
(we note that here $H_{e_1}(n)$ are the Hankel determinants (of
dimension $n$) associated with the formal power series
$P_{e_0}(z_1)$). By Algorithm 7.2.1 in \cite[p.~248]{jt80}, which
is based on Theorem 7.14 \cite[p.~244--246]{jt80} and given by
W.~Gragg \cite{gr74}, there exist numbers $p_{ne_1}$ and
$q_{ne_1},$ $n\ge1,$ such that $p_{ne_1}\ne0,$ $n\ge1,$ and
\begin{align*}
P_{e_0}(z_1)=\sum_{n=1}^\infty{c}_{ne_1}z_1^n\sim
\frac{p_{e_1}z_1}{1+q_{e_1}z_1}{\atop-}\frac{p_{2e_1}z_1^2}{1+q_{2e_1}z_1}{\atop-}\frac{p_{3e_1}z_1^2}{1+q_{3e_1}z_1}{\atop-}\ldots=F_{e_0}(z_1),
\end{align*}
where the symbol '$\sim$' means the correspondence between
$P_{e_0}(z_1)$ and $F_{e_0}(z_1)$ (at the origin). The
coefficients $p_{ne_1}$ and $q_{ne_1},$ $n\ge1,$ are given by the
formulas
\begin{align*}
p_{(n+1)e_1}=\frac{\sigma_{ne_1}}{\sigma_{(n-1)e_1}},\quad
q_{(n+1)e_1}=\tau_{(n-1)e_1}-\tau_{ne_1}\quad\mbox{for}\quad
n\ge0,
\end{align*}
where
\begin{align*}
\sigma_{ne_1}=\sum_{r=0}^nc_{(2n+1-r)e_1}B_{ne_1}(r),\quad\tau_{ne_1}\sigma_{ne_1}=\sum_{r=0}^nc_{(2n+2-r)e_1}B_{ne_1}(r),\\
\end{align*}
and for $1\le r\le n+1,$
\begin{align*}
B_{(n+1)e_1}(r)=B_{ne_1}(r)+q_{(n+1)e_1}B_{ne_1}(r-1)-p_{(n+1)e_1}B_{(n-1)e_1}(r-2)
\end{align*}
with the initial conditions
\begin{align*}
\sigma_{-e_1}&=B_{0e_1}(0)=B_{(n+1)e_1}(0)=1,\\
\tau_{-e_1}&=B_{(n-1)e_1}(-1)=B_{ne_1}(n+1)=0.
\end{align*}

Thus, we can write
\begin{align*}
L({\bf{z}})\sim
F_{e_0}(z_1)+\sum_{{i_1}=2}^Nc_{e_{i_1}}z_{i_1}R_{e_{i_1}}({\bf{z}}).
\end{align*}

Step 1.3: Let $H_{e_{i_1}}(n)\ne0$ for $2\le{{i_1}}\le{N}$ and
$n\ge1,$ where
\begin{equation}\label{eq:4}
H_{e_{i_1}}(n)=\left\vert
\begin{array}{cccc}
c_{e_{i_1}}&c_{2e_{i_1}}&\ldots&c_{ne_{i_1}}\\
c_{e_{2e_{i_1}}}&c_{3e_{i_1}}&\ldots&c_{(n+1)e_{i_1}}\\\hdotsfor{4}\\
c_{ne_{i_1}}&c_{(n+1)e_{i_1}}&\ldots&c_{(2n-1)e_{i_1}}
\end{array}\right\vert.
\end{equation}
By Algorithm 7.2.1 in \cite[p.~748]{jt80}, for each
$2\le{i_1}\le{N}$ there exist numbers $p'_{ne_{i_1}}$ and
$q'_{ne_{i_1}},$ $n\ge1,$ such that $p'_{ne_{i_1}}\ne0,$ $n\ge1,$
and
\begin{align*}
\sum_{n=1}^\infty{c}_{ne_{i_1}}z_{i_1}^n\sim\frac{p'_{e_{i_1}}z_{i_1}}{1+q'_{e_{i_1}}z_{i_1}}{\atop-}
\frac{p'_{2e_{i_1}}z_{i_1}^2}{1+q'_{2e_{i_1}}z_{i_1}}{\atop-}
\frac{p'_{3e_{i_1}}z_{i_1}^2}{1+q'_{3e_{i_1}}z_{i_1}}{\atop-}\ldots.
\end{align*}
The coefficients $p'_{ne_{i_1}}$ and $q'_{ne_{i_1}},$ $n\ge1,$ are
given by the formulas
\begin{align*}
p'_{(n+1)e_{i_1}}=\frac{\sigma_{ne_{i_1}}}{\sigma_{(n-1)e_{i_1}}},\quad
q'_{(n+1)e_{i_1}}=\tau_{(n-1)e_{i_1}}-\tau_{ne_{i_1}}\quad\mbox{for}\quad
n\ge0,
\end{align*}
where
\begin{equation}\label{eq:5}
\begin{split}
\sigma_{ne_{i_1}}&=\sum_{r=0}^nc_{(2n+1-r)e_{i_1}}B_{ne_{i_1}}(r),\\
\tau_{ne_{i_1}}\sigma_{ne_{i_1}}&=\sum_{r=0}^nc_{(2n+2-r)e_{i_1}}B_{ne_{i_1}}(r)
\end{split}
\end{equation}
and for $1\le r\le n+1,$
\begin{equation}\label{eq:6}
B_{(n+1)e_{i_1}}(r)=B_{ne_{i_1}}(r)+q_{(n+1)e_{i_1}}B_{ne_{i_1}}(r-1)-p_{(n+1)e_{i_1}}B_{(n-1)e_{i_1}}(r-2)
\end{equation}
with the initial conditions
\begin{equation}\label{eq:7}
\begin{split}
\sigma_{-e_{i_1}}&=B_{0e_{i_1}}(0)=B_{(n+1)e_{i_1}}(0)=1,\\
\tau_{-e_{i_1}}&=B_{(n-1)e_{i_1}}(-1)=B_{ne_{i_1}}(n+1)=0.
\end{split}
\end{equation}

Since for $2\le{{i_1}}\le{N},$
\[c_{e_{i_1}}=c_{e_{i_1}}B_{0e_{i_1}}(0)=\frac{\sigma_{0e_{i_1}}}{\sigma_{-e_{i_1}}}=p'_{e_{i_1}},\] we set $p_{e_{i_1}}=p'_{e_{i_1}},$
$2\le{{i_1}}\le{N}.$

Thus,
\begin{align*}
L({\bf{z}})\sim F_{e_0}(z_1)+
\sum_{i_1=2}^Np_{e_{i_1}}z_{i_1}R_{e_{i_1}}({\bf{z}}).
\end{align*}

Step 1.4: For each $2\le{i_1}\le{N}$ by
\begin{equation}\label{eq:8}
R'_{e_{i_1}}({\bf{z}})=\sum_{\substack{\vert{\bf{k}}\vert\ge0\\k_j=0,\;i_1+1\le{j}\le{N}}}c_{{\bf{k}}}^{e_{i_1}}{\bf{z}}^{{\bf{k}}}
\end{equation}
we denote a formal multiple power series reciprocal to
$R_{e_{i_1}}({\bf{z}}).$ The coefficients of the formal multiple
power series (\ref{eq:8}) are uniquely determined by the
recurrence relations
\begin{equation}\label{eq:9}
c_{{\bf{k}}}^{e_{i_1}}=-\sum_{\vert{\bf{r}}\vert=1}^{\vert{\bf{k}}\vert}c_{{\bf{k}}-{\bf{r}}}^{e_{i_1}}\frac{c_{{\bf{r}}+e_{i_1}}}{c_{e_{i_1}}},
\end{equation}
where $c_{\bf0}^{e_{i_1}}=1,$ moreover,
$c_{{\bf{k}}}^{e_{i_1}}=0,$ if there exists an index $j,$
$1\le{j}\le{N},$ such that $k_j<0.$

Thus, we can write
\begin{align*}
L({\bf{z}})\sim F_{e_0}(z_1)+
\sum_{i_1=2}^N\frac{p_{e_{i_1}}z_{i_1}}{R'_{e_{i_1}}({\bf{z}})}.
\end{align*}

The next construction of the multidimensional \emph{A}-fraction
with independent variables will be carried out using the ideas
outlined in Steps 1.1--1.4.

Step 2.1: 
Let $c_{e_{i_2}}^{e_{i_1}}\ne0$ for $2\le{i_2}\le{i_1}$ and
$2\le{i_1}\le{N}.$ In addition, for the formation of partial
denominators of the multidimensional \emph{A}-fraction with
independent variables we set the following conditions
\begin{align*}
c_{ne_{i_{2}}}^{e_{i_1}}=0\quad\mbox{for}\quad1\le{i_2}\le{i_1}-1,\;2\le{i_1}\le{N}\quad\mbox{and}\quad
n\ge1.
\end{align*}
Then for each $2\le{i_1}\le{N}$ we can rewrite the formal multiple
power series (\ref{eq:8}) in the form
\begin{align*}
R'_{e_{i_1}}({\bf{z}})=1+c_{e_{i_1}}^{e_{i_1}}z_{i_1}+z_{i_1}P_{e_{i_1}}(z_1)+
\sum_{i_2=2}^{i_1}c_{e_{i(2)}}^{e_{i_1}}z_{i_1}z_{i_2}R_{e_{i(2)}}({\bf{z}}),
\end{align*}
where
\begin{align*}
P_{e_{i_1}}(z_1)=\sum_{n=1}^\infty{c}_{e_{i_1}+ne_1}^{e_{i_1}}z_1^n,\quad
R_{e_{i(2)}}({\bf{z}})=\sum_{\substack{\vert{\bf{k}}\vert\ge0\\k_j=0,\;i_2+1\le{j}\le{N}}}
\frac{c_{{\bf{k}}+e_{i(2)}}^{e_{i_1}}}{c_{e_{i(2)}}^{e_{i_1}}}{\bf{z}}^{{\bf{k}}}.
\end{align*}
Since
\begin{align*}
c_{e_{i_1}}^{e_{i_1}}=-\frac{c_{2e_{i_1}}}{c_{e_{i_1}}}
=-\frac{c_{2e_{i_1}}B_{0e_{i_1}}(0)}{\sigma_{0e_{i_1}}}=-\tau_{0e_{i_1}}=q'_{e_{i_1}},\quad
2\le i_1\le N,
\end{align*}
we set $q_{e_{i_1}}=q'_{e_{i_1}},$ $2\le{{i_1}}\le{N}.$

Thus,
\begin{align*}
L({\bf{z}})\sim
F_{e_0}(z_1)+\sum_{i_1=2}^N\frac{p_{e_{i_1}}z_{i_1}}{1+q_{e_{i_1}}z_{i_1}+z_{i_1}P_{e_{i_1}}(z_1)}{\atop+}
\sum_{i_2=2}^{i_1}c_{e_{i(2)}}^{e_{i_1}}z_{i_1}z_{i_2}R_{e_{i(2)}}({\bf{z}}).
\end{align*}

Step 2.2: 
Let $H_{e_{i_1}+e_1}^{e_{i_1}}(n)\ne0$ for $2\le{i_1}\le{N}$ and
$n\ge1,$ where
\begin{equation*}
H_{e_{i_1}+e_1}^{e_{i_1}}(n)=\left\vert
\begin{array}{cccc}
c_{e_{i_1}+e_1}^{e_{i_1}}&c_{e_{i_1}+2e_1}^{e_{i_1}}&\ldots&c_{e_{i_1}+ne_1}^{e_{i_1}}\\
c_{e_{i_1}+2e_1}^{e_{i_1}}&c_{e_{i_1}+3e_1}^{e_{i_1}}&\ldots&c_{e_{i_1}+(n+1)e_1}^{e_{i_1}}\\\hdotsfor{4}\\
c_{e_{i_1}+ne_1}^{e_{i_1}}&c_{e_{i_1}+(n+1)e_1}^{e_{i_1}}&\ldots&c_{e_{i_1}+(2n-1)e_1}^{e_{i_1}}
\end{array}\right\vert.
\end{equation*}
By Algorithm 7.2.1 in \cite[p.~748]{jt80}, for each
$2\le{i_1}\le{N}$ there exist numbers $p_{e_{i_1}+ne_1}$ and
$q_{e_{i_1}+ne_1},$ $n\ge1,$ such that $p_{e_{i_1}+ne_1}\ne0,$
$n\ge1,$ and
\begin{align*}
\sum_{n=1}^\infty{c}_{e_{i_1}+ne_1}^{e_{i_1}}z_1^n\sim\frac{p_{e_{i_1}+e_1}z_1}{1+q_{e_{i_1}+e_1}z_1}
{\atop-}\frac{p_{e_{i_1}+2e_1}z_1^2}{1+q_{e_{i_1}+2e_1}z_1}
{\atop-}\frac{p_{e_{i_1}+3e_1}z_1^2}{1+q_{e_{i_1}+3e_1}z_1}{\atop-}\ldots=F_{e_{i_1}}(z_1).
\end{align*}
The coefficients $p_{e_{i_1}+ne_1}$ and $q_{e_{i_1}+ne_1},$
$n\ge1,$ are given by the formulas
\begin{align*}
p_{e_{i_1}+(n+1)e_1}=\frac{\sigma_{e_{i_1}+ne_1}^{e_{i_1}}}{\sigma_{e_{i_1}+(n-1)e_1}^{e_{i_1}}},\quad
q_{e_{i_1}+(n+1)e_1}=\tau_{e_{i_1}+(n-1)e_1}^{e_{i_1}}-\tau_{e_{i_1}+ne_1}^{e_{i_1}}\quad\mbox{for}\quad
n\ge0,
\end{align*}
where
\begin{align*}
\sigma_{e_{i_1}+ne_1}^{e_{i_1}}&=\sum_{r=0}^nc_{e_{i_1}+(2n+1-r)e_1}^{e_{i_1}}B_{e_{i_1}+ne_1}^{e_{i_1}}(r),\\
\tau_{e_{i_1}+ne_1}^{e_{i_1}}\sigma_{e_{i_1}+ne_1}^{e_{i_1}}&=\sum_{r=0}^nc_{e_{i_1}+(2n+2-r)e_1}^{e_{i_1}}B_{e_{i_1}+ne_1}^{e_{i_1}}(r),
\end{align*}
and for $1\le r\le n+1,$
\begin{align*}
B_{e_{i_1}+(n+1)e_1}^{e_{i_1}}(r)&=B_{e_{i_1}+ne_1}^{e_{i_1}}(r)+q_{e_{i_1}+(n+1)e_1}B_{e_{i_1}+ne_1}^{e_{i_1}}(r-1)
\\&\quad-p_{e_{i_1}+(n+1)e_1}B_{e_{i_1}+(n-1)e_1}^{e_{i_1}}(r-2)
\end{align*}
with the initial conditions
\begin{align*}
\sigma_{e_{i_1}-e_1}^{e_{i_1}}&=B_{e_{i_1}+0e_1}^{e_{i_1}}(0)=B_{e_{i_1}+(n+1)e_1}^{e_{i_1}}(0)=1,\\
\tau_{e_{i_1}-e_1}^{e_{i_1}}&=B_{e_{i_1}+(n-1)e_1}^{e_{i_1}}(-1)=B_{e_{i_1}+ne_1}^{e_{i_1}}(n+1)=0.
\end{align*}

Thus, we can write
\begin{align*}
L({\bf{z}})\sim
F_{e_0}(z_1)+\sum_{i_1=2}^N\frac{p_{e_{i_1}}z_{i_1}}{1+q_{e_{i_1}}z_{i_1}+z_{i_1}F_{e_{i_1}}(z_1)}{\atop+}
\sum_{i_2=2}^{i_1}c_{e_{i(2)}}^{e_{i_1}}z_{i_1}z_{i_2}R_{e_{i(2)}}({\bf{z}}).
\end{align*}

Step 2.3: 
Let $H_{e_{i_1}+e_{i_2}}^{e_{i_1}}(n)\ne0$ for
$2\le{i_2}\le{i_1}-1,$ $2\le{i_1}\le{N}$ and $n\ge1,$ where
\begin{align*}
H_{e_{i_1}+e_{i_2}}^{e_{i_1}}(n)=\left\vert
\begin{array}{cccc}
c_{e_{i_1}+e_{i_2}}^{e_{i_1}}&c_{e_{i_1}+2e_{i_2}}^{e_{i_1}}&\ldots&c_{e_{i_1}+ne_{i_2}}^{e_{i_1}}\\
c_{e_{i_1}+2e_{i_2}}^{e_{i_1}}&c_{e_{i_1}+3e_{i_2}}^{e_{i_1}}&\ldots&c_{e_{i_1}+(n+1)e_{i_2}}^{e_{i_1}}\\\hdotsfor{4}\\
c_{e_{i_1}+ne_{i_2}}^{e_{i_1}}&c_{e_{i_1}+(n+1)e_{i_2}}^{e_{i_1}}&\ldots&c_{e_{i_1}+(2n-1)e_{i_2}}^{e_{i_1}}
\end{array}\right\vert.
\end{align*}
Then, by Algorithm 7.2.1 in \cite[p.~748]{jt80}, for each
$2\le{i_2}\le{i_1}-1$ and $2\le{i_1}\le{N}$ there exist numbers
$p'_{e_{i_1}+ne_{i_2}}$ and $q'_{e_{i_1}+ne_{i_2}},$ $n\ge1,$ such
that $p'_{e_{i_1}+ne_{i_2}}\ne0,$ $n\ge1,$ and
\begin{align*}
\sum_{n=1}^\infty{c}_{e_{i_1}+ne_{i_2}}^{e_{i_1}}z_{i_2}^n\sim\frac{p'_{e_{i_1}+e_{i_2}}z_{i_2}}{1+q'_{e_{i_1}+e_{i_2}}z_{i_2}}{\atop-}
\frac{p'_{e_{i_1}+2e_{i_2}}z_{i_2}^2}{1+q'_{e_{i_1}+2e_{i_2}}z_{i_2}}{\atop-}
\frac{p'_{e_{i_1}+3e_{i_2}}z_{i_2}^2}{1+q'_{e_{i_1}+3e_{i_2}}z_{i_2}}{\atop-}\ldots.
\end{align*}
The coefficients $p'_{e_{i_1}+ne_{i_2}}$ and
$q'_{e_{i_1}+ne_{i_2}},$ $n\ge1,$ are given by the formulas
\begin{align*}
p'_{e_{i_1}+(n+1)e_{i_2}}=\frac{\sigma_{e_{i_1}+ne_{i_2}}^{e_{i_1}}}{\sigma_{e_{i_1}+(n-1)e_{i_2}}^{e_{i_1}}},\quad
q'_{e_{i_1}+(n+1)e_{i_2}}=\tau_{e_{i_1}+(n-1)e_{i_2}}^{e_{i_1}}-\tau_{e_{i_1}+ne_{i_2}}^{e_{i_1}}\quad\mbox{for}\quad
n\ge0,
\end{align*}
where
\begin{align*}
\sigma_{e_{i_1}+ne_{i_2}}^{e_{i_1}}=\sum_{r=0}^nc_{e_{i_1}+(2n+1-r)e_{i_2}}^{e_{i_1}}B_{e_{i_1}+ne_{i_2}}^{e_{i_1}}(r),\\
\tau_{e_{i_1}+ne_{i_2}}^{e_{i_1}}\sigma_{e_{i_1}+ne_{i_2}}^{e_{i_1}}=
\sum_{r=0}^nc_{e_{i_1}+(2n+2-r)e_{i_2}}^{e_{i_1}}B_{e_{i_1}+ne_{i_2}}^{e_{i_1}}(r),
\end{align*}
and for $1\le r\le n+1,$
\begin{align*}
B_{e_{i_1}+(n+1)e_{i_2}}^{e_{i_1}}(r)&=B_{e_{i_1}+ne_{i_2}}^{e_{i_1}}(r)+q_{e_{i_1}+(n+1)e_{i_2}}B_{e_{i_1}+ne_{i_2}}^{e_{i_1}}(r-1)
\\&\quad-p_{e_{i_1}+(n+1)e_{i_2}}B_{e_{i_1}+(n-1)e_{i_2}}^{e_{i_1}}(r-2)
\end{align*}
with the initial conditions
\begin{align*}
\sigma_{e_{i_1}-e_{i_2}}^{e_{i_1}}&=B_{e_{i_1}+0e_{i_2}}^{e_{i_1}}(0)=B_{e_{i_1}+(n+1)e_{i_2}}^{e_{i_1}}(0)=1,\\
\tau_{e_{i_1}-e_{i_2}}^{e_{i_1}}&=B_{e_{i_1}+(n-1)e_{i_2}}^{e_{i_1}}(-1)=B_{e_{i_1}+ne_{i_2}}^{e_{i_1}}(n+1)=0.
\end{align*}
Since for $2\le{i_2}\le{i_1}-1,$ $2\le{i_1}\le{N},$
\begin{align*}
c_{e_{i(2)}}^{e_{i_1}}=c_{e_{i(2)}}^{e_{i_1}}
B_{e_{i_1}+0e_{i_2}}^{e_{i_1}}(0)=\frac{\sigma_{e_{i_1}+0e_{i_2}}^{e_{i_1}}}{\sigma_{e_{i_1}-e_{i_2}}^{e_{i_1}}}=p'_{e_{i(2)}},
\end{align*}
and for $2\le{i_1}\le{N},$
\begin{align*}
c_{2e_{i_1}}^{e_{i_1}}&=-\frac{c_{e_{i_1}}^{e_{i_1}}c_{2e_{i_1}}+c_{3e_{i_2}}}{c_{e_{i_1}}}=
-\frac{c_{3e_{i_1}}c_{e_{i_1}}-(c_{2e_{i_1}})^2}{(c_{e_{i_1}})^2}=
-\frac{c_{3e_{i_1}}+c_{2e_{i_1}}q_{e_{i_1}}B_{0e_{i_1}}(0)}{c_{e_{i_1}}}\\&=
-\frac{c_{3e_{i_1}}B_{e_{i_1}}(0)+c_{2e_{i_1}}B_{e_{i_1}}(1)}{c_{e_{i_1}}B_{0e_{i_1}}(0)}=
-\frac{\sigma_{e_{i_1}}}{\sigma_{0e_{i_1}}}=-p'_{2e_{i_1}},
\end{align*}
we put
\begin{align*}
p_{e_{i(2)}}=p'_{e_{i(2)}},\quad p_{2e_{i_1}}=p'_{2e_{i_1}},\quad
2\le{i_2}\le{i_1}-1,\quad 2\le{i_1}\le{N}.
\end{align*}
Thus,
\begin{align*}
L({\bf{z}})\sim
F_{e_0}(z_1)+\sum_{i_1=2}^N\frac{p_{e_{i_1}}z_{i_1}}{1+q_{e_{i_1}}z_{i_1}+z_{i_1}F_{e_{i_1}}(z_1)}{\atop+}
\sum_{i_2=2}^{i_1}(-1)^{\delta_{i_1,i_2}}p_{e_{i(2)}}z_{i_1}z_{i_2}R_{e_{i(2)}}({\bf{z}}).
\end{align*}

Step 2.4: 
Let for each $2\le{i_2}\le{i_1}$ and $2\le{i_1}\le{N}$
\begin{gather}\label{eq:10}
R'_{e_{i(2)}}({\bf{z}})=\sum_{\substack{\vert{\bf{k}}\vert\ge0\\k_j=0,\;i_2+1\le{j}\le{N}}}c_{{\bf{k}}}^{e_{i(2)}}{\bf{z}}^{\bf{k}}
\end{gather}
be reciprocal to the formal multiple power series
$R_{e_{i(2)}}({\bf{z}}).$ It is  known that the coefficients
$c_{{\bf{k}}}^{e_{i(2)}},$ $\vert{\bf{k}}\vert\ge1,$ $k_j=0,$
$i_k+1\le{j}\le{N},$ of (\ref{eq:10}) are uniquely determined by a
recurrence formula
\begin{align*}
c_{{\bf{k}}}^{e_{i(2)}}=-\sum_{\vert{\bf{r}}\vert=1}^{\vert{\bf{k}}\vert}c_{{\bf{k}}-{\bf{r}}}^{e_{i(2)}}
\frac{c_{{\bf{r}}+e_{i(2)}}^{e_{i_1}}}{c_{e_{i(2)}}^{e_{i_1}}},
\end{align*}
where $c_{\bf{0}}^{e_{i(2)}}=1,$ moreover,
$c_{{\bf{k}}}^{e_{i(2)}}=0,$ if there exists an index $j$ such
that $1\le{j}\le{N}$ and that $k_j<0.$ Then
\begin{align*}
L({\bf{z}})\sim
F_{e_0}(z_1)+\sum_{i_1=2}^N\frac{p_{e_{i_1}}z_{i_1}}{1+q_{e_{i_1}}z_{i_1}+z_{i_1}F_{e_{i_1}}(z_1)}{\atop+}
\sum_{i_2=2}^{i_1}\frac{(-1)^{\delta_{i_1,i_2}}p_{e_{i(2)}}z_{i_1}z_{i_2}}{R'_{e_{i(2)}}({\bf{z}})}.
\end{align*}

Let us continue the construction of the multidimensional
\emph{A}-fraction with independent variables.

Step 3.1: 
Let
\begin{align*}
c_{e_{i_2}+e_{i_3}}^{e_{i(2)}}\ne0\quad\mbox{for}\quad
1\le{i_3}\le{i_2},\quad 2\le{i_2}\le{i_1},\quad 2\le{i_1}\le{N},
\end{align*}
and
\begin{align*}
c_{ne_{i_3}}^{e_{i(2)}}=0\quad\mbox{for}\quad
1\le{i_3}\le{i_2}-1,\quad 2\le{i_2}\le{i_1},\quad
2\le{i_1}\le{N},\quad\mbox{and}\quad n\ge1.
\end{align*}
Then for each $2\le{i_2}\le{i_1}$ and $2\le{i_1}\le{N}$ we have
\begin{align*}
R'_{e_{i(2)}}({\bf{z}})=1+c_{e_{i_2}}^{e_{i(2)}}z_{i_2}+z_{i_2}P_{e_{i(2)}}(z_1)+
\sum_{i_3=2}^{i_2}c_{e_{i_2,i_3}}^{e_{i(2)}}z_{i_2}z_{i_3}R_{e_{i(3)}}({\bf{z}}),
\end{align*}
where
\begin{align*}
P_{e_{i(2)}}(z_1)=\sum_{n=1}^\infty{c}_{e_{i_2}+ne_1}^{e_{i(2)}}z_1^n,\quad
{R}_{e_{i(3)}}({\bf{z}})=\sum_{\substack{\vert{\bf{k}}\vert\ge0\\{\bf{k}}_j=0,\;i_3+1\le{j}\le{N}}}
\frac{c_{{\bf{k}}+e_{i_2,i_3}}^{e_{i(2)}}}{c_{e_{i_2,i_3}}^{e_{i(2)}}}{\bf{z}}^{{\bf{k}}}.
\end{align*}
Since for $2\le{i_2}\le{i_1}-1,$ $2\le{i_1}\le{N},$
\begin{align*}
c_{e_{i_2}}^{e_{i(2)}}=-\frac{c_{e_{i_2}+e_{i(2)}}^{e_{i_1}}}{c_{e_{i(2)}}^{e_{i_1}}}
=-\frac{c_{e_{i_1}+2e_{i_2}}^{e_{i_1}}B_{e_{i_1}+0e_{i_2}}(0)}{\sigma_{e_{i_1}+0e_{i_2}}^{e_{i_1}}}=-\tau_{e_{i_1}+0e_{i_2}}^{e_{i_1}}=q'_{e_{i(2)}},
\end{align*}
and for $2\le{i_1}\le{N},$
\begin{align*}
c_{e_{i_1}}^{2e_{i_1}}&=-\frac{c_{3e_{i_1}}^{e_{i_1}}}{c_{2e_{i_1}}^{e_{i_1}}}=-\frac{c_{2e_{i_1}}^{e_{i_1}}c_{2e_{i_1}}+
c_{e_{i_1}}^{e_{i_1}}c_{3e_{i_1}}+c_{4e_{i_1}}}{c_{e_{i_1}}^{e_{i_1}}c_{2e_{i_1}}+c_{3e_{i_1}}}=\frac{c_{2e_{i_1}}}{c_{e_{i_1}}}-
\frac{c_{4e_{i_1}}c_{e_{i_1}}-c_{3e_{i_1}}c_{2e_{i_1}}}{c_{3e_{i_1}}c_{e_{i_1}}-(c_{2e_{i_1}})^2}\\&=
\frac{c_{2e_{i_1}}}{c_{e_{i_1}}}-
\frac{c_{4e_{i_1}}+c_{3e_{i_1}}q_{e_{i_1}}B_{0e_{i_1}}(0)}{c_{3e_{i_1}}B_{e_{i_1}}(0)+c_{2e_{i_1}}q_{e_{i_1}}B_{0e_{i_1}}(0)}=
\frac{c_{2e_{i_1}}}{c_{e_{i_1}}}-
\frac{c_{4e_{i_1}}+c_{3e_{i_1}}B_{e_{i_1}}(1)}{c_{3e_{i_1}}B_{e_{i_1}}(0)+c_{2e_{i_1}}B_{e_{i_1}}(1)}
\\&=
\frac{c_{2e_{i_1}}B_{0e_{i_1}}(0)}{\sigma_{0e_{i_1}}}-
\frac{c_{4e_{i_1}}B_{e_{i_1}}(0)+c_{3e_{i_1}}B_{e_{i_1}}(1)}{\sigma_{e_{i_1}}}=\tau_{0e_{i_1}}-\tau_{e_{i_1}}
=q'_{2e_{i_1}},
\end{align*}
we set
\begin{align*}
q_{e_{i(2)}}=q'_{e_{i(2)}},\quad q_{2e_{i_1}}=q'_{2e_{i_1}},\quad
2\le{i_2}\le{i_1}-1,\quad 2\le{i_1}\le{N}.
\end{align*}
Thus,
\begin{align*}
L({\bf{z}})&\sim
F_{e_0}(z_1)+\sum_{i_1=2}^N\frac{p_{e_{i_1}}z_{i_1}}{1+q_{e_{i_1}}z_{i_1}+z_{i_1}F_{e_{i_1}}(z_1)}\\&\quad{\atop+}
\sum_{i_2=2}^{i_1}\frac{(-1)^{\delta_{i_1,i_2}}p_{e_{i(2)}}z_{i_1}z_{i_2}}{1+q_{e_{i(2)}}z_{i_2}+z_{i_2}P_{e_{i(2)}}(z_1)}{\atop+}
\sum_{i_3=2}^{i_2}c_{e_{i_2,i_3}}^{e_{i(2)}}z_{i_2}z_{i_3}R_{e_{i(3)}}({\bf{z}}).
\end{align*}

Step 3.2: 
Let $H_{e_{i_2}+e_1}^{e_{i(2)}}(n)\ne0$ for $2\le{i_2}\le{i_1},$
$2\le{i_1}\le{N}$ and $n\ge1,$ where
\begin{align*}
H_{e_{i_2}+e_1}^{e_{i(2)}}(n)=\left\vert
\begin{array}{cccc}
c_{e_{i_2}+e_1}^{e_{i(2)}}&c_{e_{i_2}+2e_1}^{e_{i(2)}}&\ldots&c_{e_{i_2}+ne_1}^{e_{i(2)}}\\
c_{e_{i_2}+2e_1}^{e_{i(2)}}&c_{e_{i_2}+3e_1}^{e_{i(2)}}&\ldots&c_{e_{i_2}+(n+1)e_1}^{e_{i(2)}}\\\hdotsfor{4}\\
c_{e_{i_2}+ne_1}^{e_{i(2)}}&c_{e_{i_2}+(n+1)e_1}^{e_{i(2)}}&\ldots&c_{e_{i_2}+(2n-1)e_1}^{e_{i(2)}}
\end{array}\right\vert.
\end{align*}
Then, by Algorithm 7.2.1 in \cite[p.~748]{jt80}, for each
$2\le{i_2}\le{i_1}$ and $2\le{i_1}\le{N}$ there exist numbers
$p_{e_{i(2)}+ne_1}$ and $q_{e_{i(2)}+ne_1},$ $n\ge1,$ such that
$p_{e_{i(2)}+ne_1}\ne0$ for $n\ge1$ and
\begin{align*}
\sum_{n=1}^\infty{c}_{e_{i_2}+ne_1}^{e_{i(2)}}z_1^n\sim\frac{p_{e_{i(2)}+e_1}z_1}{1+q_{e_{i(2)}+e_1}z_1}{\atop-}
\frac{p_{e_{i(2)}+2e_1}z_1^2}{1+q_{e_{i(2)}+2e_1}z_1}{\atop-}
\frac{p_{e_{i(2)}+3e_1}z_1^2}{1+q_{e_{i(2)}+3e_1}z_1}{\atop-}\ldots=F_{e_{i(2)}}(z_1).
\end{align*}
The coefficients $p_{e_{i(2)}+ne_1}$ and $q_{e_{i(2)}+ne_1},$
$n\ge1,$ are given by the formulas for $n\ge0,$
\begin{align*}
p_{e_{i(2)}+(n+1)e_1}=\frac{\sigma_{e_{i(2)}+ne_1}^{e_{i(2)}}}{\sigma_{e_{i(2)}+(n-1)e_1}^{e_{i(2)}}},\quad
q_{e_{i(2)}+(n+1)e_1}=\tau_{e_{i(2)}+(n-1)e_1}^{e_{i(2)}}-\tau_{e_{i(2)}+ne_1}^{e_{i(2)}},
\end{align*}
where
\begin{align*}
\sigma_{e_{i(2)}+ne_1}^{e_{i(2)}}&=\sum_{r=0}^nc_{e_{i(2)}+(2n+1-r)e_1}^{e_{i(2)}}B_{e_{i(2)}+ne_1}^{e_{i(2)}}(r),\\
\tau_{e_{i(2)}+ne_1}^{e_{i(2)}}\sigma_{e_{i(2)}+ne_1}^{e_{i(2)}}&=\sum_{r=0}^nc_{e_{i(2)}+(2n+2-r)e_1}^{e_{i(2)}}B_{e_{i(2)}+ne_1}^{e_{i(2)}}(r),
\end{align*}
and for $1\le r\le n+1,$
\begin{align*}
B_{e_{i(2)}+(n+1)e_1}^{e_{i(2)}}(r)&=B_{e_{i(2)}+ne_1}^{e_{i(2)}}(r)+q_{e_{i(2)}+(n+1)e_1}B_{e_{i(2)}+ne_1}^{e_{i(2)}}(r-1)
\\&\quad-p_{e_{i(2)}+(n+1)e_1}B_{e_{i(2)}+(n-1)e_1}^{e_{i(2)}}(r-2)
\end{align*}
with the initial conditions
\begin{align*}
\sigma_{e_{i(2)}-e_1}^{e_{i(2)}}&=B_{e_{i(2)}+0e_1}^{e_{i(2)}}(0)=B_{e_{i(2)}+(n+1)e_1}^{e_{i(2)}}(0)=1,\\
\tau_{e_{i(2)}-e_1}^{e_{i(2)}}&=B_{e_{i(2)}+(n-1)e_1}^{e_{i(2)}}(-1)=B_{e_{i(2)}+ne_1}^{e_{i(2)}}(n+1)=0.
\end{align*}

Thus,
\begin{align*}
L({\bf{z}})&\sim
F_{e_0}(z_1)+\sum_{i_1=2}^N\frac{p_{e_{i_1}}z_{i_1}}{1+q_{e_{i_1}}z_{i_1}+z_{i_1}F_{e_{i_1}}(z_1)}\\&\quad{\atop+}
\sum_{i_2=2}^{i_1}\frac{(-1)^{\delta_{i_1,i_2}}p_{e_{i(2)}}z_{i_1}z_{i_2}}{1+q_{e_{i(2)}}z_{i_2}+z_{i_2}F_{e_{i(2)}}(z_1)}{\atop+}
\sum_{i_3=2}^{i_2}c_{e_{i_2,i_3}}^{e_{i(2)}}z_{i_2}z_{i_3}R_{e_{i(3)}}({\bf{z}}).
\end{align*}

Step 3.3: 
Let $H_{e_{i_2}+e_{i_3}}^{e_{i(2)}}(n)\ne0$ for
$2\le{i_3}\le{i_2}-1,$ $2\le{i_2}\le{i_1},$ $2\le{i_1}\le{N}$ and
$n\ge1,$ where
\begin{align*}
H_{e_{i_2}+e_{i_3}}^{e_{i(2)}}(n)=\left\vert
\begin{array}{cccc}
c_{e_{i_2}+e_{i_3}}^{e_{i(2)}}&c_{e_{i_2}+2e_{i_3}}^{e_{i(2)}}&\ldots&c_{e_{i_2}+ne_{i_3}}^{e_{i(2)}}\\
c_{e_{i_2}+2e_{i_3}}^{e_{i(2)}}&c_{e_{i_2}+3e_{i_3}}^{e_{i(2)}}&\ldots&c_{e_{i_2}+(n+1)e_{i_3}}^{e_{i(2)}}\\\hdotsfor{4}\\
c_{e_{i_2}+ne_{i_3}}^{e_{i(2)}}&c_{e_{i_2}+(n+1)e_{i_3}}^{e_{i(2)}}&\ldots&c_{e_{i_2}+(2n-1)e_{i_3}}^{e_{i(2)}}
\end{array}\right\vert.
\end{align*}
Then, by Algorithm 7.2.1 in \cite[p.~748]{jt80}, for each
$2\le{i_3}\le{i_2}-1,$ $2\le{i_2}\le{i_1}$ and $2\le{i_1}\le{N}$
there exist numbers $p'_{e_{i(2)}+ne_{i_3}}$ and
$q'_{e_{i(2)}+ne_{i_3}},$ $n\ge1,$ such that
$p'_{e_{i(2)}+ne_{i_3}}\ne0,$ $n\ge1,$ and
\begin{align*}
\sum_{n=1}^\infty{c}_{e_{i_2}+ne_{i_3}}^{e_{i(2)}}z_{i_3}^n\sim\frac{p'_{e_{i(3)}}z_{i_3}}{1+q'_{e_{i(3)}}z_{i_3}}{\atop-}
\frac{p'_{e_{i(2)}+2e_{i_3}}z_{i_3}^2}{1+q'_{e_{i(2)}+2e_{i_3}}z_{i_3}}{\atop-}
\frac{p'_{e_{i(2)}+3e_{i_3}}z_{i_3}^2}{1+q'_{e_{i(2)}+3e_{i_3}}z_{i_3}}{\atop-}\ldots.
\end{align*}
The coefficients $p'_{e_{i(2)}+ne_{i_3}}$ and
$q'_{e_{i(2)}+ne_{i_3}},$ $n\ge1,$ are given by the formulas for
$n\ge0,$
\begin{align*}
p'_{e_{i(2)}+(n+1)e_{i_3}}=\frac{\sigma_{e_{i(2)}+ne_{i_3}}^{e_{i(2)}}}{\sigma_{e_{i(2)}+(n-1)e_{i_3}}^{e_{i(2)}}},\quad
q'_{e_{i(2)}+(n+1)e_{i_3}}=\tau_{e_{i(2)}+(n-1)e_{i_3}}^{e_{i(2)}}-\tau_{e_{i(2)}+ne_{i_3}}^{e_{i(2)}},
\end{align*}
where
\begin{align*}
\sigma_{e_{i(2)}+ne_{i_3}}^{e_{i(2)}}&=\sum_{r=0}^nc_{e_{i(2)}+(2n+1-r)e_{i_3}}^{e_{i(2)}}B_{e_{i(2)}+ne_{i_3}}^{e_{i(2)}}(r),\\
\tau_{e_{i(2)}+ne_{i_3}}^{e_{i(2)}}\sigma_{e_{i(2)}+ne_{i_3}}^{e_{i(2)}}&=
\sum_{r=0}^nc_{e_{i(2)}+(2n+2-r)e_{i_3}}^{e_{i(2)}}B_{e_{i(2)}+ne_{i_3}}^{e_{i(2)}}(r),
\end{align*}
and for $1\le r\le n+1,$
\begin{align*}
B_{e_{i(2)}+(n+1)e_{i_3}}^{e_{i(2)}}(r)&=B_{e_{i(2)}+ne_{i_3}}^{e_{i(2)}}(r)+q_{e_{i(2)}+(n+1)e_{i_3}}B_{e_{i(2)}+ne_{i_3}}^{e_{i(2)}}(r-1)
\\&\quad-p_{e_{i(2)}+(n+1)e_{i_3}}B_{e_{i(2)}+(n-1)e_{i_3}}^{e_{i(2)}}(r-2)
\end{align*}
with the initial conditions
\begin{align*}
\sigma_{e_{i(2)}-e_{i_3}}^{e_{i(2)}}&=B_{e_{i(2)}+0e_{i_3}}^{e_{i(2)}}(0)=B_{e_{i(2)}+(n+1)e_{i_3}}^{e_{i(2)}}(0)=1,\\
\tau_{e_{i(2)}-e_{i_3}}^{e_{i(2)}}&=B_{e_{i(2)}+(n-1)e_{i_3}}^{e_{i(2)}}(-1)=B_{e_{i(2)}+ne_{i_3}}^{e_{i(2)}}(n+1)=0.
\end{align*}
Since for $2\le{i_3}\le{i_2}-1,$ $2\le{i_2}\le{i_1},$
$2\le{i_1}\le{N},$
\begin{align*}
c_{e_{i_2}+e_{i_3}}^{e_{i(2)}}=c_{e_{i_2}+e_{i_3}}^{e_{i(2)}}B_{e_{i(2)}+0e_{i_3}}^{e_{i(2)}}(0)
=\frac{\sigma_{e_{i(2)}+0e_{i_3}}^{e_{i(2)}}}{\sigma_{e_{i(2)}-e_{i_3}}^{e_{i(2)}}}=p'_{e_{i(3)}},
\end{align*}
and for $2\le{i_2}\le{i_1}-1,$ $3\le{i_1}\le{N},$
\begin{align*}
c_{2e_{i_2}}^{e_{i(2)}}&=-c_{e_{i_2}}^{e_{i(2)}}\frac{c_{e_{i_1}+2e_{i_2}}^{e_{i_1}}}{c_{e_{i(2)}}^{e_{i_1}}}-
\frac{c_{e_{i_1}+3e_{i_2}}^{e_{i_1}}}{c_{e_{i(2)}}^{e_{i_1}}}\\&=-\frac{(c_{e_{i_1}+2e_{i_2}}^{e_{i_1}})^2-
c_{e_{i(2)}}^{e_{i_1}}c_{e_{i_1}+3e_{i_2}}^{e_{i_1}}}{(c_{e_{i(2)}}^{e_{i_1}})^2}
=-\frac{c_{e_{i_1}+3e_{i_2}}^{e_{i_1}}+c_{e_{i_1}+2e_{i_2}}^{e_{i_1}}q_{e_{i(2)}}}
{c_{e_{i(2)}}^{e_{i_1}}}
\\&=-\frac{c_{e_{i_1}+3e_{i_2}}^{e_{i_1}}B_{e_{i(2)}}^{e_{i_1}}(0)+c_{e_{i_1}+2e_{i_2}}^{e_{i_1}}B_{e_{i(2)}}^{e_{i_1}}(1)}
{c_{e_{i(2)}}^{e_{i_1}}B_{e_{i_1}+0e_{i_2}}^{e_{i_1}}(0)}
=-\frac{\sigma_{e_{i(2)}}^{e_{i_1}}}{\sigma_{e_{i_1}+0e_{i_2}}^{e_{i_1}}}
=-p'_{e_{i_1}+2e_{i_2}},
\end{align*}
(note that the coefficient $c_{2e_{i_2}}^{e_{i(2)}}$ is possible
only if $N\ge3$ and, of course, the appearance of this coefficient
here and similar others in the following steps depends on the
number $N$), and for $2\le{i_1}\le{N},$
\begin{align*}
c_{2e_{i_1}}^{2e_{i_1}}&=\frac{(c_{3e_{i_1}}^{e_{i_1}})^2-c_{2e_{i_1}}^{e_{i_1}}c_{4e_{i_1}}^{e_{i_1}}}{(c_{2e_{i_1}}^{e_{i_1}})^2}\\&=
-\frac{c_{e_{i_1}}(c_{e_{i_1}}c_{3e_{i_1}}c_{5e_{i_1}}+2c_{2e_{i_1}}c_{3e_{i_1}}c_{4e_{i_1}}-c_{3e_{i_1}}^3-
c_{e_{i_1}}c_{4e_{i_1}}^2-c_{2e_{i_1}}^2c_{5e_{i_1}})}{(c_{e_{i_1}}c_{3e_{i_1}}-c_{2e_{i_1}}^2)^2}\\&
=-\frac{c_{5e_{i_1}}+c_{4e_{i_1}}(q_{e_{i_1}}+q_{2e_{i_1}})+c_{3e_{i_1}}(q_{e_{i_1}}q_{2e_{i_1}}-p_{2e_{i_1}})}
{c_{3e_{i_1}}+c_{2e_{i_1}}q_{e_{i_1}}B_{0e_{i_1}}(0)}\\&
=-\frac{c_{5e_{i_1}}B_{2e_{i_1}}(0)+c_{4e_{i_1}}B_{2e_{i_1}}(1)+c_{3e_{i_1}}B_{2e_{i_1}}(2)}
{c_{3e_{i_1}}B_{e_{i_1}}(0)+c_{2e_{i_1}}B_{e_{i_1}}(1)}=
-\frac{\sigma_{2e_{i_1}}}{\sigma_{e_{i_1}}}=-p'_{3e_{i_1}},
\end{align*}
we put
\begin{align*}
p_{e_{i(3)}}&=p'_{e_{i(3)}},\quad 2\le{i_3}\le{i_2}-1,\quad 2\le{i_2}\le{i_1},\quad 2\le{i_1}\le{N},\\
p_{e_{i_1}+2e_{i_2}}&=p'_{e_{i_1}+2e_{i_2}},\quad 2\le{i_3}\le{i_2}-1,\quad 3\le{i_1}\le{N},\\
p_{3e_{i_1}}&=p'_{3e_{i_1}},\quad 2\le{i_1}\le{N}.
\end{align*}

Thus,
\begin{align*}
L({\bf{z}})&\sim
F_{e_0}(z_1)+\sum_{i_1=2}^N\frac{p_{e_{i_1}}z_{i_1}}{1+q_{e_{i_1}}z_{i_1}+z_{i_1}F_{e_{i_1}}(z_1)}\\&\quad{\atop+}
\sum_{i_2=2}^{i_1}\frac{(-1)^{\delta_{i_1,i_2}}p_{e_{i(2)}}z_{i_1}z_{i_2}}{1+q_{e_{i(2)}}z_{i_2}+z_{i_2}F_{e_{i(2)}}(z_1)}{\atop+}
\sum_{i_3=2}^{i_2}(-1)^{\delta_{i_2,i_3}}p_{e_{i(3)}}z_{i_2}z_{i_3}R_{e_{i(3)}}({\bf{z}}).
\end{align*}

Step 3.4: 
We obtain
\begin{align*}
L({\bf{z}})&\sim
F_{e_0}(z_1)+\sum_{i_1=2}^N\frac{p_{e_{i_1}}z_{i_1}}{1+q_{e_{i_1}}z_{i_1}+z_{i_1}F_{e_{i_1}}(z_1)}\\&\quad{\atop+}
\sum_{i_2=2}^{i_1}\frac{(-1)^{\delta_{i_1,i_2}}p_{e_{i(2)}}z_{i_1}z_{i_2}}{1+q_{e_{i(2)}}z_{i_2}+z_{i_2}F_{e_{i(2)}}(z_1)}{\atop+}
\sum_{i_3=2}^{i_2}\frac{(-1)^{\delta_{i_2,i_3}}p_{e_{i(3)}}z_{i_2}z_{i_3}}{R'_{e_{i(3)}}({\bf{z}})},
\end{align*}
where for each $2\le{i_3}\le{i_2},$ $2\le{i_2}\le{i_1},$ and
$2\le{i_1}\le{N},$
\begin{equation}\label{eq:11}
R'_{e_{i(3)}}({\bf{z}})=\sum_{\substack{\vert{\bf{k}}\vert\ge0\\
{\bf{k}}_j=0,\;i_3+1\le{j}\le{N}}}c_{{\bf{k}}}^{e_{i(3)}}{\bf{z}}^{\bf{k}}
\end{equation}
is reciprocal to the $R_{e_{i(3)}}({\bf{z}}).$ The coefficients of
(\ref{eq:11}) calculated as follows
\begin{align*}
c_{{\bf{k}}}^{e_{i(3)}}=-\sum_{\vert{\bf{r}}\vert=1}^{\vert{\bf{k}}\vert}c_{{\bf{k}}-{\bf{r}}}^{e_{i(3)}}
\frac{c_{{\bf{r}}+e_{i_2,i_3}}^{e_{i(2)}}}{c_{e_{i_2,i_3}}^{e_{i(2)}}},
\end{align*}
where $c_{\bf{0}}^{e_{i(3)}}=1,$ moreover,
$c_{{\bf{k}}}^{e_{i(3)}}=0,$ if there exists an index $j$ such
that $1\le{j}\le{N}$ and $k_j<0.$

Further construction of the multidimensional \emph{A}-fraction
with independent variables (\ref{eq:2}) consists of gradually
applying steps similar to Steps 2.1--2.4 to all formal multiple
power series in the denominators of the ending partial quotients
of the finite branches of the branched continued fraction.

 As a result, computing the coefficients
$c_{{\bf{k}}}^{e_{i_1}},$ $\vert{\bf{k}}\vert\ge1,$
$k_j=0,\;i_1+1\le{j}\le{N},$ $2\le i_1\le N,$ by the recurrence
formula (\ref{eq:9}) and the coefficients
$c_{{\bf{k}}}^{e_{i(k)}},$ $\vert{\bf{k}}\vert\ge1,$
$k_j=0,\;i_k+1\le{j}\le{N},$ $k\ge2,$ $2\le i_p\le i_{p-1},$ $1\le
p\le k,$ by the recurrence formula
\begin{equation}\label{eq:12}
c_{{\bf{k}}}^{e_{i(k)}}=-\sum_{\vert{\bf{r}}\vert=1}^{\vert{\bf{k}}\vert}c_{{\bf{k}}-{\bf{r}}}^{e_{i(k)}}
\frac{c_{{\bf{r}}+e_{i_{k-1},i_k}}^{e_{i(k-1)}}}{c_{e_{i_{k-1},i_k}}^{e_{i(k-1)}}},
\end{equation}
where $c_{\bf{0}}^{e_{i(k)}}=1,$ moreover,
$c_{{\bf{k}}}^{e_{i(k)}}=0,$ if there exists an index  $j,$
$1\le{j}\le{N},$ such that $k_j<0,$ provided that for $1\le i_1\le
N$ and $n\ge1,$
\begin{equation}\label{eq:13}
H_{e_{i_1}}(n)\ne0,
\end{equation}
where $H_{e_{i_1}}(n)$ defined by (\ref{eq:4}), and provided that
for each $1\le i_{k+1}\le i_k-1,$ $2\le i_p\le i_{p-1},$ $1\le
p\le k,$ $k\ge1$ and $n\ge1,$
\begin{equation}\label{eq:14}
c_{ne_{i_{k+1}}}^{e_{i(k)}}=0,\quad
H_{e_{i_k,i_{k+1}}}^{e_{i(k)}}(n)\ne0,
\end{equation}
where $H_{e_{i_k,i_{k+1}}}^{e_{i(k)}}(n)$ defined by
\begin{align*}
H_{e_{i_k,i_{k+1}}}^{e_{i(k)}}(n)=\left\vert
\begin{array}{cccc}
c_{e_{i_k,i_{k+1}}}^{e_{i(k)}}&c_{e_{i_k}+2e_{i_{k+1}}}^{e_{i(k)}}&\ldots&c_{e_{i_k}+ne_{i_{k+1}}}^{e_{i(k)}}\\
c_{e_{i_k}+2e_{i_{k+1}}}^{e_{i(k)}}&c_{e_{i_k}+3e_{i_{k+1}}}^{e_{i(k)}}&\ldots&c_{e_{i_k}+(n+1)e_{i_{k+1}}}^{e_{i(k)}}\\\hdotsfor{4}\\
c_{e_{i_k}+ne_{i_{k+1}}}^{e_{i(k)}}&c_{e_{i_k}+(n+1)e_{i_{k+1}}}^{e_{i(k)}}&\ldots&c_{e_{i_k}+(2n-1)e_{i_{k+1}}}^{e_{i(k)}}
\end{array}\right\vert,
\end{align*}
for the formal multiple power series (\ref{eq:3}) we obtain
multidimensional \emph{A}-fraction with independent variables
(\ref{eq:2}), where the $p_{e_{i(k)}}$ and $q_{e_{i(k)}}$ for all
$e_{i(k)}\in\mathcal{E}_k,$ $k\ge1,$ defined by the following
formulas:
\begin{equation}
\begin{split}\label{eq:15}
p_{(n+1)e_{i_1}}&=\frac{\sigma_{ne_{i_1}}}{\sigma_{(n-1)e_{i_1}}},\\
q_{(n+1)e_{i_1}}&=\tau_{(n-1)e_{i_1}}-\tau_{ne_{i_1}},
\end{split}
\end{equation}
where $1\le i_1\le N,$ $n\ge0,$ and $\sigma_{ne_{i_1}},$
$\tau_{ne_{i_1}},$ $n\ge-1,$ defined by
(\ref{eq:5})--(\ref{eq:7}),
\begin{equation}
\begin{split}\label{eq:16}
p_{e_{i(k)}+(n+1)e_{i_{k+1}}}&=\frac{\sigma_{e_{i(k)}+ne_{i_{k+1}}}^{e_{i(k)}}}{\sigma_{e_{i(k)}+(n-1)e_{i_{k+1}}}^{e_{i(k)}}},\\
q_{e_{i(k)}+(n+1)e_{i_{k+1}}}&=\tau_{e_{i(k)}+(n-1)e_{i_{k+1}}}^{e_{i(k)}}-\tau_{e_{i(k)}+ne_{i_{k+1}}}^{e_{i(k)}},
\end{split}
\end{equation}
where $2\le i_p\le i_{p-1},$ $1\le p\le k,$ $1\le i_{k+1}\le
i_k-1,$ $k\ge1,$ $n\ge0,$
\begin{align*}
\sigma_{e_{i(k)}+ne_{i_{k+1}}}^{e_{i(k)}}&=\sum_{r=0}^nc_{e_{i(k)}+(2n+1-r)e_{i_{k+1}}}^{e_{i(k)}}B_{e_{i(k)}+ne_{i_{k+1}}}^{e_{i(k)}}(r),\\
\tau_{e_{i(k)}+ne_{i_{k+1}}}^{e_{i(k)}}\sigma_{e_{i(k)}+ne_{i_{k+1}}}^{e_{i(k)}}&=
\sum_{r=0}^nc_{e_{i(k)}+(2n+2-r)e_{i_{k+1}}}^{e_{i(k)}}B_{e_{i(k)}+ne_{i_{k+1}}}^{e_{i(k)}}(r),
\end{align*}
and for $1\le r\le n+1,$
\begin{align*}
B_{e_{i(k)}+(n+1)e_{i_{k+1}}}^{e_{i(k)}}(r)&=B_{e_{i(k)}+ne_{i_{k+1}}}^{e_{i(k)}}(r)+q_{e_{i(k)}+(n+1)e_{i_{k+1}}}B_{e_{i(k)}+ne_{i_{k+1}}}^{e_{i(k)}}(r-1)
\\&\quad-p_{e_{i(k)}+(n+1)e_{i_{k+1}}}B_{e_{i(k)}+(n-1)e_{i_{k+1}}}^{e_{i(k)}}(r-2)
\end{align*}
with the initial conditions
\begin{align*}
\sigma_{e_{i(2)}-e_{i_3}}^{e_{i(2)}}&=B_{e_{i(2)}+0e_{i_3}}^{e_{i(2)}}(0)=B_{e_{i(2)}+(n+1)e_{i_3}}^{e_{i(2)}}(0)=1,\\
\tau_{e_{i(2)}-e_{i_3}}^{e_{i(2)}}&=B_{e_{i(2)}+(n-1)e_{i_3}}^{e_{i(2)}}(-1)=B_{e_{i(2)}+ne_{i_3}}^{e_{i(2)}}(n+1)=0.
\end{align*}

Thus, we have constructed the recurrence algorithm for computing
the coefficients of the multidimensional \emph{A}-fraction with
independent variables (\ref{eq:2}) in terms of the given formal
multiple power series (\ref{eq:3}).

Further, we show that the constructed multidimensional
\emph{A}-fraction with independent variables (\ref{eq:2})
corresponds at ${\bf{z}}={\bf{0}}$ to the formal multiple power
series (\ref{eq:3}).

Let $\{f_n({\bf{z}})\}$ be a sequence of approximants of the
multidimensional \emph{A}-fraction with independent variables
(\ref{eq:2}). Using formulas (\ref{eq:9}), (\ref{eq:12}),
(\ref{eq:15}), and (\ref{eq:16}) we curtail the $f_n({\bf{z}})$
for $n\ge1.$

Note that according to the described above algorithm  for $e_0$
and for all $e_{i(k)}$ such that $2\le i_p\le i_{p-1},$ $1\le p\le
k$ and $k\ge1$ the continued fraction
\begin{align*}
F_{e_{i(k)}}(z_1)=\frac{p_{e_{i(k)}+e_1}z_1}{1+q_{e_{i(k)}+e_1}z_1}{\atop-}
\frac{p_{e_{i(k)} +2e_1}z_1^2}{1+q_{e_{i(k)}+2e_1}z_1}{\atop-}
\frac{p_{e_{i(k)}
+3e_1}z_1^2}{1+q_{e_{i(k)}+3e_1}z_1}{\atop-}\ldots
\end{align*}
corresponds at the origin to the formal power series
\begin{align*}
P_{e_{i(k)}}(z_1)=\sum_{r=1}^\infty
c_{e_{i_k}+re_1}^{e_{i(k)}}z_1^r
\end{align*}
and the order of correspondence is $\nu_n=2n+1.$ It follows that
for $e_0$ and for each $e_{i(k)}$ such that $2\le i_p\le i_{p-1},$
$1\le p\le k,$ $k\ge1$ and for $n\ge2$ the finite continued
fraction
\begin{align*}
F_{e_{i(k)}}^{(n)}(z_1)=\frac{p_{e_{i(k)}+e_1}z_1}{1+p_{e_{i(k)}+e_1}z_1}{\atop-}
\frac{p_{e_{i(k)}+2e_1}z_1^2}{1+q_{e_{i(k)}+2e_1}z_1}{\atop-}\ldots{\atop-}
\frac{p_{e_{i(k)}+ne_1}z_1^2}{1+q_{e_{i(k)}+ne_1}z_1}
\end{align*}
has formal power series expansion
\begin{align*}
P_{e_{i(k)}}^{(n)}(z_1)=\sum_{r=1}^{2n}{c}_{e_{i_k}+re_1}^{e_{i(k)}}z_1^r+O(z_1^{2n+1}),
\end{align*}
where $c_{re_1}^{e_0}=c_{re_1}$ for $1\le r\le 2n$ and $n\ge1,$
$O(z_1^p)$ is a symbolic mark for some formal power series, whose
minimal degree of terms is not less than $p,$ $p\ge3.$

Now for $n=1$ we have
\begin{align*}
f_1({\bf{z}})&=\sum_{i_1=1}^N\frac{p_{e_{i_1}}z_{i_1}}{1+q_{e_{i_1}}z_{i_1}}=
\sum_{i_1=1}^N\frac{c_{e_{i_1}}z_{i_1}}{1+c_{e_{i_1}}^{e_{i_1}}z_{i_1}}
=\sum_{i_1=1}^N\frac{c_{e_{i_1}}z_{i_1}}{1}{\atop+}\sum_{i_2=1}^{i_1}c_{e_{i_2}}^{e_{i_1}}z_{i_2}\\
&=\sum_{i_1=1}^Nc_{e_{i_1}}z_{i_1}\bigg(1+\sum_{i_2=1}^{i_1}\frac{c_{e_{i(2)}}}{c_{e_{i_1}}}z_{i_2}+O({\bf{z}}^2)\bigg)\\
&=\sum_{i_1=1}^Nc_{e_{i_1}}z_{i_1}+\sum_{i_1=1}^Nz_{i_1}\bigg(\sum_{i_2=1}^{i_1}c_{e_{i(2)}}z_{i_2}\bigg)+O({\bf{z}}^3)
=\sum_{\vert{\bf{k}}\vert=1}^2c_{\bf{k}}{\bf{z}}^{\bf{k}}+O({\bf{z}}^3),
\end{align*}
where $O({\bf{z}}^p)$ is a symbolic mark for some formal multiple
power series, whose minimal degree of homogeneous terms is not
less than $p,$ $p\ge2.$ Since
\begin{align*}
\sum_{\vert{\bf{k}}\vert=1}^2c_{\bf{k}}{\bf{z}}^{\bf{k}}+O({\bf{z}}^{3})-\sum_{\vert{\bf{k}}\vert\ge1}c_{\bf{k}}{\bf{z}}^{\bf{k}}=O'({\bf{z}}^{3}),
\end{align*}
where $O'({\bf{z}}^p)$ is a symbolic mark for some formal multiple
power series, whose minimal degree of homogeneous terms is not
less than $p,$ $p\ge3,$ then $f_1({\bf{z}})\sim L({\bf{z}})$ and
the order of correspondence is $\nu_1=3.$

For $n=2$ we can write
\begin{align*}
f_2({\bf{z}})&=F_{e_0}^{(2)}(z_1)+\sum_{i_1=2}^N\frac{p_{e_{i_1}}z_{i_1}}{1+q_{e_{i_1}}z_{i_1}}{\atop+}
\sum_{i_2=1}^{i_1}\frac{(-1)^{\delta_{i_1,i_2}}p_{e_{i(2)}}z_{i_1}z_{i_2}}{1+q_{e_{i(2)}}z_{i_2}}\\
&=P_{e_0}^{(2)}(z_1)+\sum_{i_1=2}^N\frac{c_{e_{i_1}}z_{i_1}}{1+c_{e_{i_1}}^{e_{i_1}}z_{i_1}}{\atop+}
\sum_{i_2=1}^{i_1}\frac{c_{e_{i(2)}}^{e_{i_1}}z_{i_1}z_{i_2}}{1+c_{e_{i_2}}^{e_{i(2)}}z_{i_2}}\\
&=P_{e_0}^{(2)}(z_1)+\sum_{i_1=2}^N\frac{c_{e_{i_1}}z_{i_1}}{1+c_{e_{i_1}}^{e_{i_1}}z_{i_1}}{\atop+}
\sum_{i_2=1}^{i_1}\frac{c_{e_{i(2)}}^{e_{i_1}}z_{i_1}z_{i_2}}{1}{\atop+}\sum_{i_3=1}^{i_2}c_{e_{i_3}}^{e_{i(2)}}z_{i_3}\\
&=P_{e_0}^{(2)}(z_1)+\sum_{i_1=2}^N\frac{c_{e_{i_1}}z_{i_1}}{1+c_{e_{i_1}}^{e_{i_1}}z_{i_1}}
{\atop+}\sum_{i_2=1}^{i_1}c_{e_{i(2)}}^{e_{i_1}}z_{i_1}z_{i_2}
\Bigg(1+\sum_{i_3=1}^{i_2}\frac{c_{e_{i(2)}+e_{i_3}}^{e_{i_1}}}{c_{e_{i(2)}}^{e_{i_1}}}z_{i_3}+O({\bf{z}}^2)\Bigg)\\
&=P_{e_0}^{(2)}(z_1)+\sum_{i_1=2}^N\frac{c_{e_{i_1}}z_{i_1}}{1}{\atop+}
\bigg(\sum_{i_2=1}^{i_1}c_{e_{i_2}}^{e_{i_1}}z_{i_2}\\
&\quad+\sum_{i_2=1}^{i_1}z_{i_2}\bigg(\sum_{i_3=1}^{i_2}c_{e_{i_2,i_3}}^{e_{i_1}}z_{i_3}\bigg)+
\sum_{i_2=1}^{i_1}z_{i_2}\bigg(\sum_{i_3=1}^{i_2}z_{i_3}\bigg(\sum_{i_4=1}^{i_3}c_{e_{i(4)}-e_{i_1}}^{e_{i_1}}z_{i_4}\bigg)\bigg)
+O({\bf{z}}^3)\bigg)\\
&=P_{e_0}^{(2)}(z_1)+\sum_{i_1=2}^Nc_{e_{i_1}}z_{i_1}\bigg(1+\sum_{i_2=1}^{i_1}\frac{c_{e_{i(2)}}}{c_{e_{i_1}}}z_{i_2}
+\sum_{i_2=1}^{i_1}z_{i_2}\bigg(\sum_{i_3=1}^{i_2}\frac{c_{e_{i(3)}}}{c_{e_{i_1}}}z_{i_3}\bigg)\\
&\quad+\sum_{i_2=1}^{i_1}z_{i_2}\bigg(\sum_{i_3=1}^{i_2}z_{i_3}\bigg(\sum_{i_4=1}^{i_3}\frac{c_{e_{i(4)}}}{c_{e_{i_1}}}z_{i_4}\bigg)\bigg)+O({\bf{z}}^4)\bigg)
=\sum_{\vert{\bf{k}}\vert=1}^4c_{\bf{k}}{\bf{z}}^{\bf{k}}+O({\bf{z}}^5).
\end{align*}
Since
\begin{align*}
\sum_{\vert{\bf{k}}\vert=1}^4c_{\bf{k}}{\bf{z}}^{\bf{k}}+O({\bf{z}}^5)-\sum_{\vert{\bf{k}}\vert\ge1}c_{\bf{k}}{\bf{z}}^{\bf{k}}=O'({\bf{z}}^{5}),
\end{align*}
then $f_2({\bf{z}})\sim L({\bf{z}})$ and $\nu_2=5.$

Next, let $n\ge3$ be an arbitrary natural number. Then we get
\begin{align*}
f_{n}({\bf{z}})&=F_{e_0}^{(n)}(z_1)+\sum_{i_1=2}^N\frac{p_{e_{i_1}}z_{i_1}}{1+q_{e_{i_1}}z_{i_1}+z_{i_1}F_{e_{i_1}}^{(n-1)}(z_1)}\\
&\quad{\atop+}\sum_{i_2=2}^{i_1}\frac{(-1)^{\delta_{i_1,i_2}}p_{e_{i(2)}}z_{i_1}z_{i_2}}{1+q_{e_{i(2)}}z_{i_2}+z_{i_2}F_{e_{i(2)}}^{(n-2)}(z_1)}
{\atop+}\ldots{\atop+}
\sum_{i_n=1}^{i_{n-1}}\dfrac{(-1)^{\delta_{i_{n-1},i_n}}p_{e_{i(n)}}z_{i_{n-1}}z_{i_n}}{1+q_{e_{i(n)}}z_{i_n}}\\&=
P_{e_0}^{(n)}(z_1)+\sum_{i_1=2}^N\frac{c_{e_{i_1}}z_{i_1}}{1+c_{e_{i_1}}^{e_{i_1}}z_{i_1}+z_{i_1}P_{e_{i-1}}^{(n-1)}(z_1)}{\atop+}
\sum_{i_2=2}^{i_1}\frac{c_{e_{i(2)}}^{e_{i_1}}z_{i_1}z_{i_2}}{1+c_{e_{i_2}}^{e_{i(2)}}z_{i_2}+z_{i_2}P_{e_{i(2)}}^{(n-2)}(z_1)}\\
&\quad{\atop+}\ldots{\atop+}
\sum_{i_{n-2}=2}^{i_{n-3}}\frac{c_{e_{i_{n-3},i_{n-2}}}^{e_{i(n-3)}}z_{i_{n-3}}z_{i_{n-2}}}{1+c_{e_{i_{n-2}}}^{e_{i(n-2)}}z_{i_{n-2}}+
z_{i_{n-2}}P_{e_{i(n-2)}}^{(2)}(z_1)}\\&\quad{\atop+}
\sum_{i_{n-1}=2}^{i_{n-2}}\dfrac{c_{e_{i_{n-2},i_{n-1}}}^{e_{i(n-2)}}z_{i_{n-2}}z_{i_{n-1}}}{1+c_{e_{i_{n-1}}}^{e_{i(n-1)}}z_{i_{n-1}}}{\atop+}
\sum_{i_n=1}^{i_{n-1}}\dfrac{c_{e_{i_{n-1},i_n}}^{e_{i(n-1)}}z_{i_{n-1}}z_{i_n}}{1+c_{e_{i_{n}}}^{e_{i(n)}}z_{i_n}}\\&=
P_{e_0}^{(n)}(z_1)+\sum_{i_1=2}^N\frac{c_{e_{i_1}}z_{i_1}}{1+c_{e_{i_1}}^{e_{i_1}}z_{i_1}+z_{i_1}P_{e_{i_1}}^{(n-1)}(z_1)}{\atop+}
\sum_{i_2=2}^{i_1}\frac{c_{e_{i(2)}}^{e_{i_1}}z_{i_1}z_{i_2}}{1+c_{e_{i_2}}^{e_{i(2)}}z_{i_2}+z_{i_2}P_{e_{i(2)}}^{(n-2)}(z_1)}\\
&\quad{\atop+}\ldots{\atop+}
\sum_{i_{n-2}=2}^{i_{n-3}}\frac{c_{e_{i_{n-3},i_{n-2}}}^{e_{i(n-3)}}z_{i_{n-3}}z_{i_{n-2}}}{1+c_{e_{i_{n-2}}}^{e_{i(n-2)}}z_{i_{n-2}}+
z_{i_{n-2}}P_{e_{i(n-2)}}^{(2)}(z_1)}\\&\quad{\atop+}
\sum_{i_{n-1}=2}^{i_{n-2}}\dfrac{c_{e_{i_{n-2},i_{n-1}}}^{e_{i(n-2)}}z_{i_{n-2}}z_{i_{n-1}}}{1+c_{e_{i_{n-1}}}^{e_{i(n-1)}}z_{i_{n-1}}}{\atop+}
\sum_{i_n=1}^{i_{n-1}}\dfrac{c_{e_{i_{n-1},i_n}}^{e_{i(n-1)}}z_{i_{n-1}}z_{i_n}}{1}{\atop+}\sum_{i_{n+1}=1}^{i_n}c_{e_{i_{n+1}}}^{e_{i(n)}}z_{i_{n+1}}
\\&
=P_{e_0}^{(n)}(z_1)+\sum_{i_1=2}^N\frac{c_{e_{i_1}}z_{i_1}}{1+c_{e_{i_1}}^{e_{i_1}}z_{i_1}+z_{i_1}P_{e_{i_1}}^{(n-1)}(z_1)}{\atop+}
\sum_{i_2=2}^{i_1}\frac{c_{e_{i(2)}}^{e_{i_1}}z_{i_1}z_{i_2}}{1+c_{e_{i_2}}^{e_{i(2)}}z_{i_2}+z_{i_2}P_{e_{i(2)}}^{(n-2)}(z_1)}\\
&\quad{\atop+}\ldots{\atop+}
\sum_{i_{n-2}=2}^{i_{n-3}}\frac{c_{e_{i_{n-3},i_{n-2}}}^{e_{i(n-3)}}z_{i_{n-3}}z_{i_{n-2}}}{1+c_{e_{i_{n-2}}}^{e_{i(n-2)}}z_{i_{n-2}}+
z_{i_{n-2}}P_{e_{i(n-2)}}^{(2)}(z_1)}{\atop+}
\sum_{i_{n-1}=2}^{i_{n-2}}\dfrac{c_{e_{i_{n-2},i_{n-1}}}^{e_{i(n-2)}}z_{i_{n-2}}z_{i_{n-1}}}{1+c_{e_{i_{n-1}}}^{e_{i(n-1)}}z_{i_{n-1}}}\\
&\quad{\atop+}\sum_{i_n=1}^{i_{n-1}}c_{e_{i_{n-1},i_n}}^{e_{i(n-1)}}z_{i_{n-1}}z_{i_n}
\Bigg(1+\sum_{i_{n+1}=1}^{i_n}\frac{c_{e_{i_{n-1},i_n,i_{n+1}}}^{e_{i(n-1)}}}{c_{e_{i_{n-1},i_n}}^{e_{i(n-1)}}}z_{i_{n+1}}+O({\bf{z}}^2)\Bigg)\end{align*}
\begin{align*}&=
P_{e_0}^{(n)}(z_1)+\sum_{i_1=2}^N\frac{c_{e_{i_1}}z_{i_1}}{1+c_{e_{i_1}}^{e_{i_1}}z_{i_1}+z_{i_1}P_{e_{i_1}}^{(n-1)}(z_1)}{\atop+}
\sum_{i_2=2}^{i_1}\frac{c_{e_{i(2)}}^{e_{i_1}}z_{i_1}z_{i_2}}{1+c_{e_{i_2}}^{e_{i(2)}}z_{i_2}+z_{i_2}P_{e_{i(2)}}^{(n-2)}(z_1)}\\
&\quad{\atop+}\ldots{\atop+}
\sum_{i_{n-2}=2}^{i_{n-3}}\frac{c_{e_{i_{n-3},i_{n-2}}}^{e_{i(n-3)}}z_{i_{n-3}}z_{i_{n-2}}}{1+c_{e_{i_{n-2}}}^{e_{i(n-2)}}z_{i_{n-2}}+
z_{i_{n-2}}P_{e_{i(n-2)}}^{(2)}(z_1)}{\atop+}\sum_{i_{n-1}=2}^{i_{n-2}}\dfrac{c_{e_{i_{n-2},i_{n-1}}}^{e_{i(n-2)}}z_{i_{n-2}}z_{i_{n-1}}}{1}\\&
\quad{\atop+}\bigg(\sum_{i_n=1}^{i_{n-1}}c_{e_{i_n}}^{e_{i(n-1)}}z_{i_n}+
\sum_{i_n=1}^{i_{n-1}}z_{i_n}\bigg(\sum_{i_{n+1}=1}^{i_n}c_{e_{i_n,i_{n+1}}}^{e_{i(n-1)}}z_{i_{n+1}}\bigg)\\&
\quad+\sum_{i_n=1}^{i_{n-1}}z_{i_n}\bigg(\sum_{i_{n+1}=1}^{i_n}z_{i_{n+1}}
\bigg(\sum_{i_{n+2}=1}^{i_{n+1}}c_{e_{i_n,i_{n+1},i_{n+2}}}^{e_{i(n-1)}}z_{i_{n+2}}\bigg)\bigg)+O({\bf{z}}^3)\bigg).
\end{align*}
Continuing this process on the final step we obtain
\begin{align*}
f_{n}({\bf{z}})&=P_{e_0}^{(n)}(z_1)+\sum_{i_1=2}^N\frac{c_{e_{i_1}}z_{i_1}}{1}{\atop+}
\bigg(\sum_{i_2=1}^{i_1}c_{e_{i_2}}^{e_{i_1}}z_{i_2}+
\sum_{i_2=1}^{i_1}z_{i_2}\bigg(\sum_{i_3=1}^{i_2}c_{e_{i_2,i_3}}^{e_{i_1}}z_{i_3}\bigg)\\&\quad+\ldots+\sum_{i_2=1}^{i_1}z_{i_2}
\bigg(\sum_{i_3=1}^{i_2}z_{i_3}\bigg(\ldots\bigg(\sum_{i_{2n-1}=1}^{i_{2n-2}}z_{i_{2n-1}}
\bigg(\sum_{i_{2n}=1}^{i_{2n-1}}c_{e_{i(2n)}-e_{i_1}}^{e_{i_1}}z_{i_{2n}}\bigg)\bigg)\ldots\bigg)\bigg)\\&\quad
+O({\bf{z}}^{2n-1})\bigg).
\end{align*}
From this we have
\begin{align*}
f_{n}({\bf{z}})&=P_{e_0}^{(n)}(z_1)+
\sum_{i_1=2}^Nc_{e_{i_1}}z_{i_1}\bigg(1+\sum_{i_2=1}^{i_1}\frac{c_{e_{i(2)}}}{c_{e_{i_1}}}z_{i_2}+
\sum_{i_2=1}^{i_1}z_{i_2}\bigg(\sum_{i_3=1}^{i_2}\frac{c_{e_{i(3)}}}{c_{e_{i_1}}}z_{i_3}\bigg)\\&\quad+\ldots+\sum_{i_2=1}^{i_1}z_{i_2}
\bigg(\ldots\bigg(\sum_{i_{2n-1}=1}^{i_{2n-2}}z_{i_{2n-1}}
\bigg(\sum_{i_{2n}=1}^{i_{2n-1}}\frac{c_{e_{i(2n)}}}{c_{e_{i_1}}}z_{i_{2n}}\bigg)\bigg)\ldots\bigg)+O({\bf{z}}^{2n})\bigg)\\&=
\sum_{\vert{\bf{k}}\vert=1}^{2n}c_{\bf{k}}{\bf{z}}^{\bf{k}}+O({\bf{z}}^{2n+1}).
\end{align*}
Since
\begin{align*}
\sum_{\vert{\bf{k}}\vert=1}^{2n}c_{\bf{k}}{\bf{z}}^{\bf{k}}+O({\bf{z}}^{2n+1})-
\sum_{\vert{\bf{k}}\vert\ge1}c_{\bf{k}}{\bf{z}}^{\bf{k}}=O'({\bf{z}}^{2n+1}),
\end{align*}
$f_n({\bf{z}})\sim L({\bf{z}})$ and $\nu_n=2n+1.$

Al last, from arbitrariness of $n$ it follows that
$f_{n}({\bf{z}})\sim{L}({\bf{z}})$ for all $n\ge1$ and that the
order of correspondence is $\nu_{n}=2n+1.$ It follows, that
$\Lambda(f_n)$ and $L({\bf{z}})$ agree for all homogeneous terms
up to and including degree $2n.$ Since
\begin{align*}
\lim_{n\to+\infty}\nu_n=\lim_{n\to+\infty}2n+1=+\infty,
\end{align*}
the multidimensional \emph{A}-fraction with independent variables
(\ref{eq:2}) corresponds at ${\bf{z}}={\bf{0}}$ to the formal
multiple power series (\ref{eq:3}).

Thus, the following theorem is true.
\begin{theorem}\label{th:1}
The multidimensional  A-fraction with independent variables
(\ref{eq:2}) corresponds at ${\bf{z}}={\bf{0}}$ to the given
formal multiple power series (\ref{eq:3}) if and only if the
conditions (\ref{eq:13}) for $1\le i_1\le N,$ $n\ge1,$ and the
conditions (\ref{eq:14}) for $1\le i_{k+1}\le i_k-1,$ $2\le i_p\le
i_{p-1},$ $1\le p\le k,$ $k\ge1,$ $n\ge1$ are satisfied.
\end{theorem}

It follows from Theorem 1 and Theorem 2 in \cite{bodm19} that the
conditions for the existence of the generalization of Gragg's
algorithm are the same as for the algorithm in \cite{bodm19}.
However, this algorithm provides a more convenient numerical
procedure for computing the coefficients of multidimensional
\emph{A}-fractions with independent variables corresponding to a
given formal multiple power series.

Finally, we detail the correspondence of multidimensional
\emph{J}-fractions with independent variables, closely related to
multidimensional \emph{A}-fractions with independent variables
(see, for example, \cite{dm21cmfth}).

In a multidimensional \text{\emph{A}-fraction} with independent
variables (\ref{eq:2}), we set $z_i=1/w_i,$ $1\le{i}\le{N},$ and
perform the equivalence transformation (see \cite[pp.
29--33]{bod86}) by setting $\rho_{e_{i(k)}}=w_{i_k},$
${e_{i(k)}}\in\mathcal{E}_k,$ $k\ge1.$ As a result, we arrive at
the multidimensional \emph{J}-fraction with independent variables
\begin{equation}\label{eq:17}
\sum_{i_1=1}^N\frac{p_{e_{i(1)}}}{q_{e_{i(1)}}+w_{i_1}}{\atop+}
\sum_{i_2=1}^{i_1}\frac{(-1)^{\delta_{i_1,i_2}}p_{e_{i(2)}}}{q_{e_{i(2)}}+w_{i_2}}{\atop+}
\sum_{i_3=1}^{i_2}\frac{(-1)^{\delta_{i_2,i_3}}p_{e_{i(3)}}}{q_{e_{i(3)}}+w_{i_3}}{\atop+}\ldots,
\end{equation}
where $p_{e_{i(k)}},$ $q_{e_{i(k)}},$ $e_{i(k)}\in\mathcal{E}_k,$
$k\ge1,$ are complex numbers and, in addition, $p_{e_{i(k)}}\ne0,$
$e_{i(k)}\in\mathcal{E}_k,$ $k\ge1.$

A sequence of rational functions $\{R_n({\bf{w}})\}$ is said to
correspond at ${\bf{w}}=\pmb{\infty}$ to a formal multiple Laurent
series
\begin{equation}\label{eq:18}
L^*({\bf{w}})=\sum_{\vert{\bf{k}}\vert\ge0}\frac{c_{\bf{k}}}{{\bf{w}}^{\bf{k}}},
\end{equation}
where $c_{\bf{k}}\in\mathbb{C},$ ${\bf{k}}\ge1,$ if the sequence
$\{R_n(1/z_1,1/z_2,\ldots,1/z_N)\}$ corresponds to a formal
multiple power series at ${\bf{z}}={\bf{0}}$ obtained from
(\ref{eq:18}) by replacing $w_i$ with $1/z_i,$ $1\le{i}\le{N}.$

A formal multiple Laurent series (\ref{eq:18}) is said to be an
asymptotic expansion of a function $R({\bf{w}})$ at
${\bf{w}}=\pmb{\infty},$ with respect to a region $D$ in
$\mathbb{C}^N,$ if for every $n\ge0$ there exist $\rho_n>0$ and
$\eta_n>0$ such that
\[\left\vert R({\bf{w}})-\sum_{\vert{\bf{k}}\vert=0}^n\frac{c_{\bf{k}}}{{\bf{w}}^{\bf{k}}}
\right\vert\le\eta_n\left(\sum_{k=1}^N\frac{1}{\vert{w_k}\vert}\right)^{n+1},\quad
\vert{w_k}\vert>\rho_n,\quad1\le k\le N,\quad{\bf{w}}\in D.\]

We denote this by
\[R({\bf{w}})\approx\sum_{\vert{\bf{k}}\vert\ge0}\frac{c_{\bf{k}}}{{\bf{w}}^{\bf{k}}},\quad w_k\to\infty,\quad1\le k\le N.\]


A multidimensional \emph{J}-fraction with independent variables
(\ref{eq:17}) is said to correspond at ${\bf{w}}=\pmb{\infty}$ to
the formal multiple Laurent series (\ref{eq:18}) if its sequence
of approximants $\{f_n^*({\bf{w}})\}$ corresponds to
$L^*({\bf{w}})$ at ${\bf{w}}=\pmb{\infty}.$

The following theorem summarizes the connections between
multidimensional \emph{A}- and \emph{J}-fractions with independent
variables
(see \cite[Teorem~3]{bodm19}).
\begin{theorem}
\label{th:2} Let $f_n({\bf{z}})$ ($f_n^*({\bf{w}})$) denote the
$n$th approximants, respectively, of the multidimensional
A-fraction with independent variables (\ref{eq:2})
(multidimensional J-fraction with independent variables
(\ref{eq:17})), where $z_i=1/w_i,$ $1\le{i}\le{N}.$ Further, let
the multidimensional A-fraction with independent variables
(\ref{eq:2}) corresponds to the formal multiple power series
(\ref{eq:3}) at ${\bf{z}}={\bf{0}}.$ Then the following assertions
are true:
\begin{enumerate}
    \item[(A)] for any natural $n,$ the equality $f_n({\bf{z}})=f_n^*({\bf{w}})$ is
    true;
    \item[(B)] the formal expansion of the $n$th approximant $f_n^*({\bf{w}})$ in the
multiple Laurent series at ${\bf{w}}=\pmb{\infty}$ has the form
\begin{align*}
f_n^*({\bf{w}})=\sum_{\vert{\bf{k}}\vert=1}^{2n}\frac{c_{{\bf{k}}}}{{\bf{w}}^{{\bf{k}}}}+
\sum_{\vert{\bf{k}}\vert\ge2n+1}\frac{c_{{\bf{k}}}^{(n)}}{{\bf{w}}^{{\bf{k}}}},\quad{n}\ge1,
\end{align*}
where $c_{{\bf{k}}}^{(n)}\in\mathbb{C},$
$\vert{\bf{k}}\vert\ge{2n+1},$ and hence, the multidimensional
J-fraction with independent variables (\ref{eq:17}) corresponds at
${\bf{w}}=\pmb{\infty}$ to the formal multiple Laurent series
\begin{align*}
L^*({\bf{w}})=\sum_{\vert{\bf{k}}\vert\ge1}\frac{c_{\bf{k}}}{{\bf{w}}^{\bf{k}}}.
\end{align*}
\end{enumerate}
\end{theorem}

It follows from Theorem~2 that the generalization of Gragg's
algorithm can also be used for computing the coefficients of
multidimensional \emph{J}-fractions with independent variables
corresponding to a given formal multiple Laurent series.

It should be noted that for multidimensional \emph{A}- and
\emph{J}-fractions with independent variables, convergence
criteria can be found in \cite{bar14,bb21,bodm19,dm17,dm21cmfth},
and convergence rate estimates in
\cite{antdm2020a,antdm2020b,bar14,bb20}. Nevertheless, the problem
of improving and developing new methods of studying the
convergence of these branched continued fractions remains open.

\section{Numerical experiments}
In this section, we will give some applications of above
constructed algorithm.

The function of two variables
\begin{equation}\label{eq:19}
F(z_1,z_2)=\arctan(z_1)+\arctan\left(\frac{z_2}{1+z_2\arctan(z_1)}\right)
\end{equation}
has a formal double power series at origin given by
\begin{align}
L(z_1,z_2)&=\sum_{k=1}^\infty\frac{(-1)^{k+1}}{2k-1}z_1^{2k-1}\notag\\
&\quad+\sum_{r=1}^\infty\frac{(-1)^{r+1}}{2r-1}\left(\sum_{s=0}^\infty
z_2^{s+1}\left(\sum_{k=1}^\infty\frac{(-1)^k}{2k-1}z_1^{2k-1}\right)^s\right)^{2r-1}.\label{eq:20}
\end{align}
Applying recurrence algorithm constructed in Section~\ref{sec:3},
we obtain the following.

Step 1.1: We have
\begin{align*}
L(z_1,z_2)&=\sum_{k=1}^\infty\frac{(-1)^{k+1}}{2k-1}z_1^{2k-1}+z_2\left(1-z_1z_2-\frac{1}{3}z_2^2+\frac{1}{3}z_1^3z_2+z_1^2z_2^2+
z_1z_2^3+\frac{1}{5}z_2^4\right.\\
&\quad\left.-\frac{1}{5}z_1^5z_2-\frac{2}{3}z_1^4z_2^2-\frac{4}{3}z_1^3z_2^3-2z_1^2z_2^4-z_1z_2^5-\frac{1}{7}z_2^6+\ldots\right).
\end{align*}

Steps 1.2 and 1.3: By formulas (\ref{eq:15}) we obtain (see also
Table~\ref{tab:1})
\begin{align*}
p_{1,0}&=p_{0,1}=1,\\
p_{k,0}&=p_{0,k}=-\frac{(k-1)^2}{(2k-3)(2k-1)},\quad k\ge2,\\
q_{k,0}&=q_{0,k}=0,\quad k\ge0.
\end{align*}
\begin{table}[h]
\begin{center}
\begin{minipage}{245pt}
\caption{Results of algorithm applied to (\ref{eq:19}) on Steps~1.2~and~1.3 for $i_1=1,2$}\label{tab:1}%
\begin{tabular}{@{}llllllll@{}}
\hline
$n$ & $p_{ne_{i_1}}$  & $q_{ne_{i_1}}$ & $\sigma_{ne_{i_1}}$ & $\tau_{ne_{i_1}}$ & $B_{ne_{i_1}}(0)$ & $B_{ne_{i_1}}(1)$ & $B_{ne_{i_1}}(2)$\\
\hline
$-1$ & & & 1 & 0\\
0    & & & 1 & 0 & 1 \\
1    & 1 & 0 & $-1/3$ & 0 & 1 & 0\\
2    & $-1/3$ & 0 &  4/45 & 1/15 & 1 & 1/3 & 1/3 \\
3    & $-4/15$ & 0 &  \\
\hline
\end{tabular}
\end{minipage}
\end{center}
\end{table}

Thus,
\begin{align*}
L(z_1,z_2)&\sim
F_0(z_1)+z_2\left(1-z_1z_2-\frac{1}{3}z_2^2+\frac{1}{3}z_1^3z_2+z_1^2z_2^2+
z_1z_2^3+\frac{1}{5}z_2^4\right.\\
&\quad\left.-\frac{1}{5}z_1^5z_2-\frac{2}{3}z_1^4z_2^2-\frac{4}{3}z_1^3z_2^3-2z_1^2z_2^4-z_1z_2^5-\frac{1}{7}z_2^6+\ldots\right),
\end{align*}
where
\begin{equation*}
F_0(z_1)=\dfrac{p_{1,0}z_1}{1}{\atop-}\dfrac{p_{2,0}z_1^2}{1}{\atop-}\dfrac{p_{3,0}z_1^2}{1}{\atop-}\ldots.
\end{equation*}

Step 1.4: By the recurrence formula (\ref{eq:9}) we obtain
\begin{align*}
L(z_1,z_2)&\sim F_0(z_1)+\frac{z_2}{R'_{0,1}(z_1,z_2)},
\end{align*}
where
\begin{align*}
R'_{0,1}(z_1,z_2)&=1+z_1z_2+\frac{1}{3}z_2^2-\frac{1}{3}z_1^3z_2-\frac{1}{3}z_1z_2^3-\frac{4}{45}z_2^4\\&\quad+
\frac{1}{5}z_1^5z_2+\frac{1}{9}z_1^3z_2^3+\frac{1}{3}z_1^2z_2^4+\frac{12}{45}z_1z_2^5+
\frac{44}{945}z_2^6+\ldots.
\end{align*}

Step 2.1: We have
\begin{align*}
L(z_1&,z_2)\sim
F_0(z_1)+\frac{z_2}{\displaystyle1+z_2\sum_{k=1}^\infty\frac{(-1)^{k+1}}{2k-1}z_1^{2k-1}+
\frac{z_2^2}{3}R_{0,2}(z_1,z_2)},
\end{align*}
where
\begin{align*}
R_{0,2}(z_1,z_2)=1-z_1z_2-\frac{4}{15}z_2^2+\frac{1}{3}z_1^3z_2+z_1^2z_2^2+\frac{12}{15}z_1z_2^3+\frac{44}{315}z_2^4+\ldots.
\end{align*}

Steps 2.2 and 2.3: By formulas (\ref{eq:16}) we obtain (see also
Table~\ref{tab:1})
\begin{align*}
p_{1,1}&=1,\\
p_{n,1}&=-\frac{(n-1)^2}{(2n-3)(2n-1)},\quad n\ge2,\\
q_{n,1}&=0,\quad n\ge0.
\end{align*}
\begin{table}[h]
\begin{center}
\begin{minipage}{245pt}
\caption{Results of algorithm applied to (\ref{eq:19}) on Steps~2.2~and~2.3}\label{tab:1}%
\begin{tabular}{@{}llllllll@{}}
\hline $n$ & $p_{n,1}$  & $q_{n,1}$ & $\sigma_{n,1}$ &
$\tau_{n,1}$ & $B_{n,1}(0)$ & $B_{n,1}(1)$ & $B_{n,1}(2)$\\
\hline
$-1$ & & & 1 & 0\\
0    & & & 1 & 0 & 1 \\
1    & 1 & 0 & $-1/3$ & 0 & 1 & 0\\
2    & $-1/3$ & 0 &  4/45 & 1/15 & 1 & 1/3 & 1/3 \\
3    & $-4/15$ & 0 &  \\
\hline
\end{tabular}
\end{minipage}
\end{center}
\end{table}

Thus,
\begin{align*}
L(z_1,z_2)\sim
F_0(z_1)+\dfrac{z_2}{\displaystyle1+z_2F_1(z_1)+\dfrac{z_2^2}{3}R_{0,2}(z_1,z_2)},
\end{align*}
where
\begin{equation*}
F_1(z_1)=\dfrac{p_{1,1}z_1}{1}{\atop-}\dfrac{p_{2,1}z_1^2}{1}{\atop-}\dfrac{p_{3,1}z_1^2}{1}{\atop-}\ldots.
\end{equation*}

Step 2.4: By the recurrence formula (\ref{eq:12}) we obtain
\begin{align*}
L(z_1,z_2)&\sim
F_0(z_1)+\dfrac{z_2}{1+z_2F_1(z_1)+\dfrac{z_2^2/3}{R'_{0,2}(z_1,z_2)}},
\end{align*}
where
\begin{align*}
R'_{0,2}(z_1,z_2)=1+z_1z_2+\frac{4}{15}z_2^2-
\frac{1}{3}z_1^3z_2-\frac{4}{15}z_1z_2^3-\frac{44}{315}z_2^4+\ldots.
\end{align*}

And so on, at the end we will get the corresponding
two-dimensional \emph{A}-fraction with independent variables of
the form
\begin{align}\label{eq:21}
F_0(z_1)+\dfrac{z_2}{1+z_2F_1(z_1)+\dfrac{-p_{0,2}z_2^2}{1+z_2F_2(z_1)+\dfrac{-p_{0,3}z_2^2}{1+{}_{\ddots}}}},
\end{align}
where for $k\ge0$
\begin{align*}
F_k(z_1)&=\dfrac{z_1}{1}{\atop-}\dfrac{p_{2,k}z_1^2}{1}{\atop-}\dfrac{p_{3,k}z_1^2}{1}{\atop-}\ldots,\\
p_{n,k}&=p_{0,n}=-\frac{(n-1)^2}{(2n-3)(2n-1)},\quad n\ge2.
\end{align*}

\begin{figure}[h]%
\centering
\begin{subfigure}[h]{.37\linewidth}
\centering
\includegraphics[width=\linewidth]{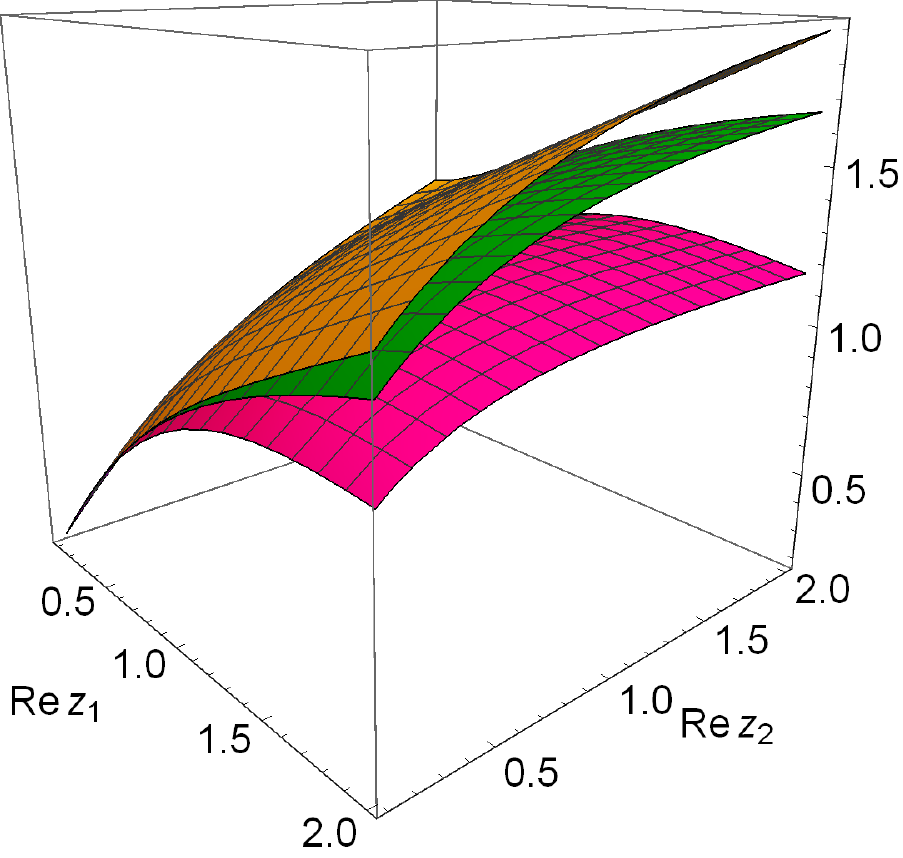}
\caption{\textcolor{rd}{$\blacksquare$} -- 2nd,
 \textcolor{gr}{$\blacksquare$}  -- $\eqref{eq:19},$
\textcolor{br}{$\blacksquare$} -- 3rd}
\end{subfigure}\quad\quad
\begin{subfigure}[h]{.37\linewidth}
\centering
\includegraphics[width=\linewidth]{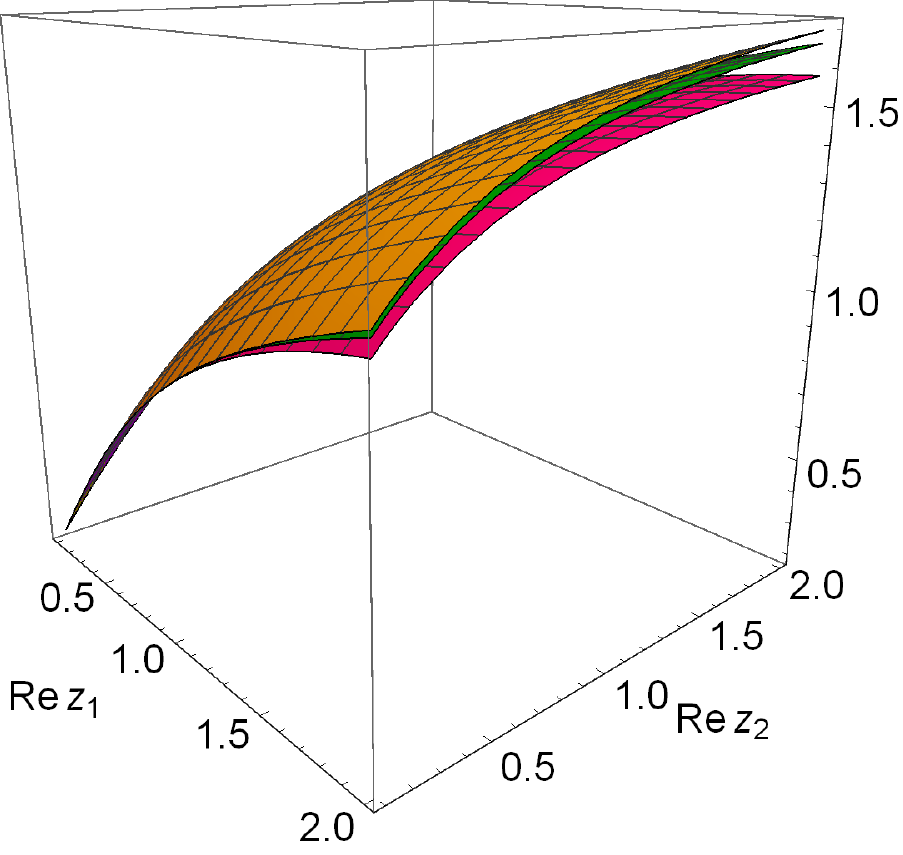}
\caption{\textcolor{rd}{$\blacksquare$} -- 4th,
 \textcolor{gr}{$\blacksquare$}  -- $\eqref{eq:19},$
\textcolor{br}{$\blacksquare$} -- 5th}
\end{subfigure}
\caption{The plots of values of the $n$th approximants of
\eqref{eq:21}}\label{fig:1}
\end{figure}

\begin{figure}[h]%
\centering
\begin{subfigure}[h]{.45\linewidth}
\centering
\includegraphics[width=\linewidth]{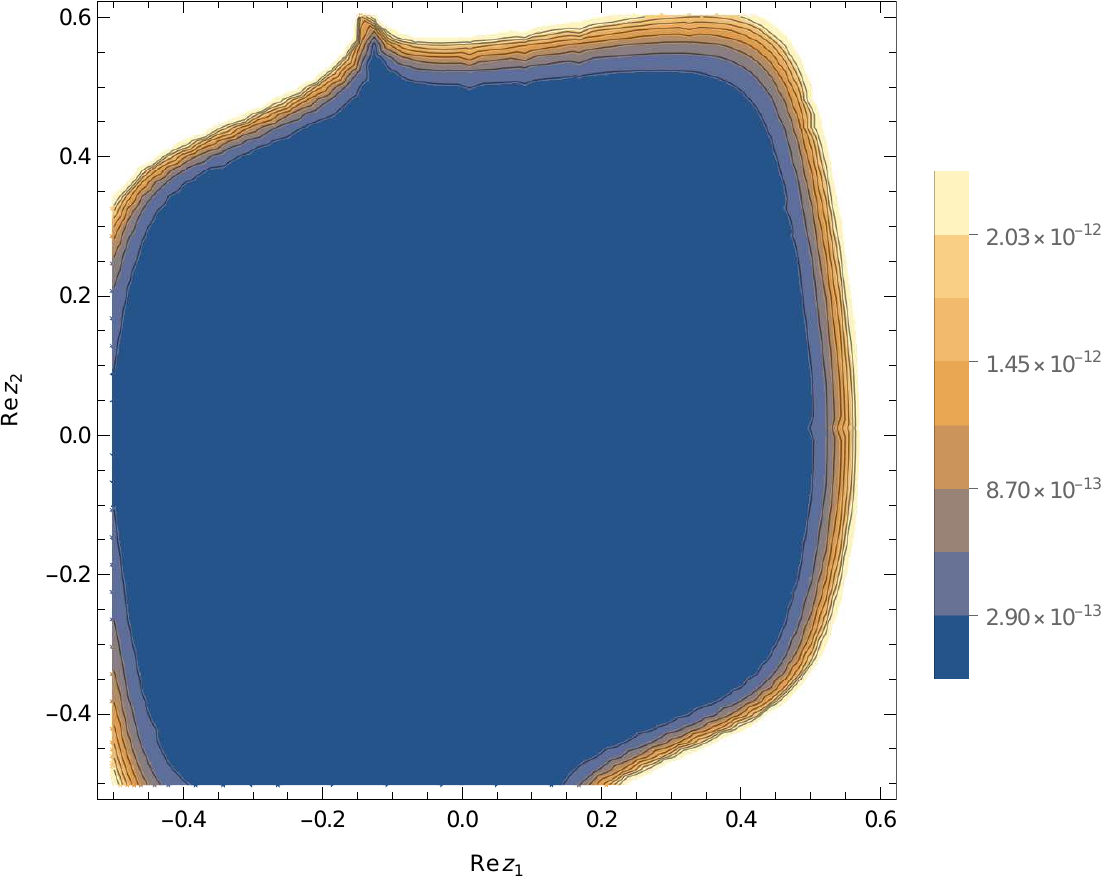}
\caption{}
\end{subfigure}
\quad\quad
\begin{subfigure}[h]{.45\linewidth}
\centering
\includegraphics[width=\linewidth]{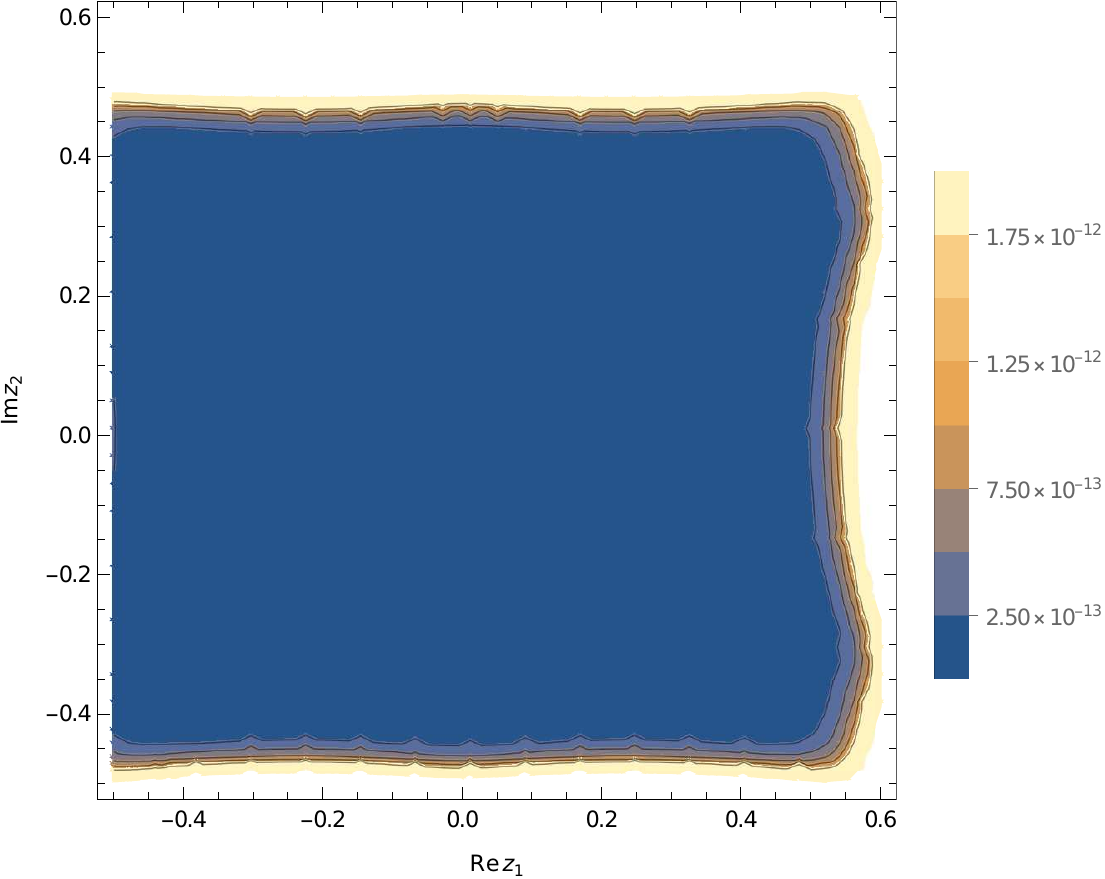}
\caption{}
\end{subfigure}\\[2ex]
\begin{subfigure}[h]{.45\linewidth}
\centering
\includegraphics[width=\linewidth]{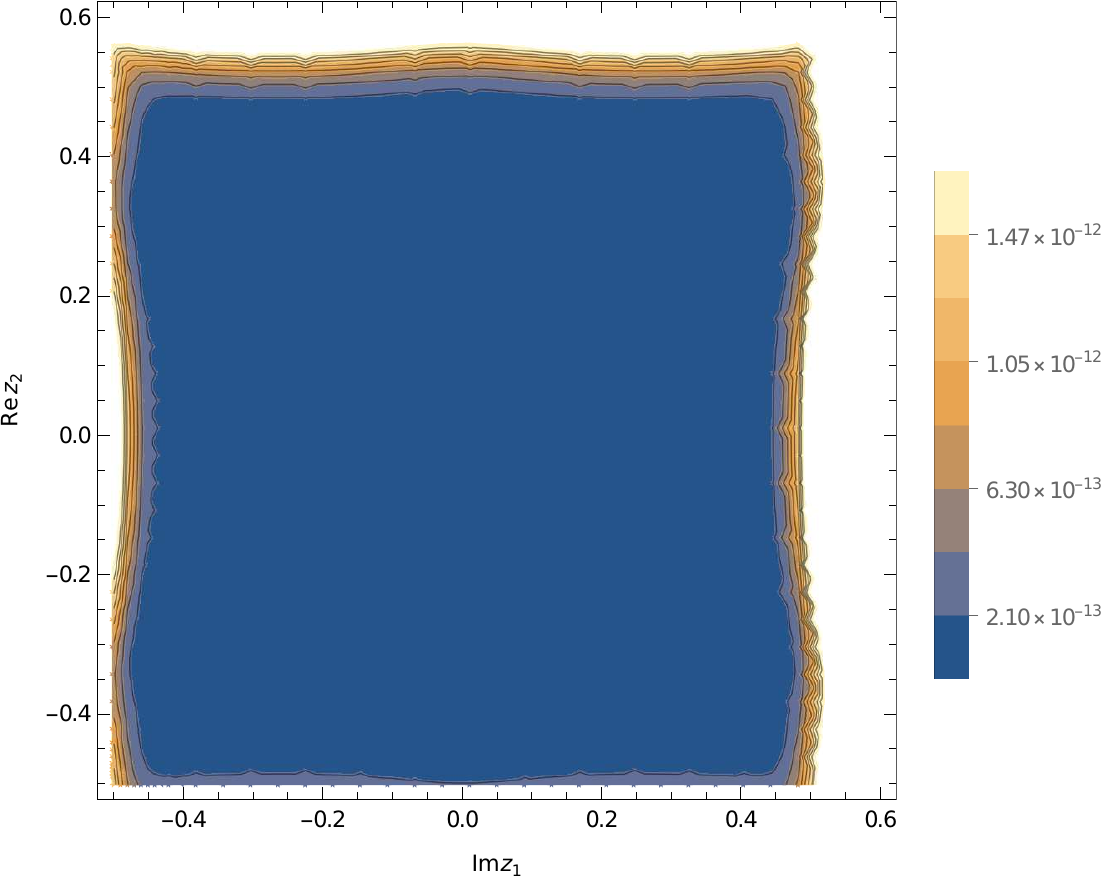}
\caption{}
\end{subfigure}\quad\quad
\begin{subfigure}[h]{.45\linewidth}
\centering
\includegraphics[width=\linewidth]{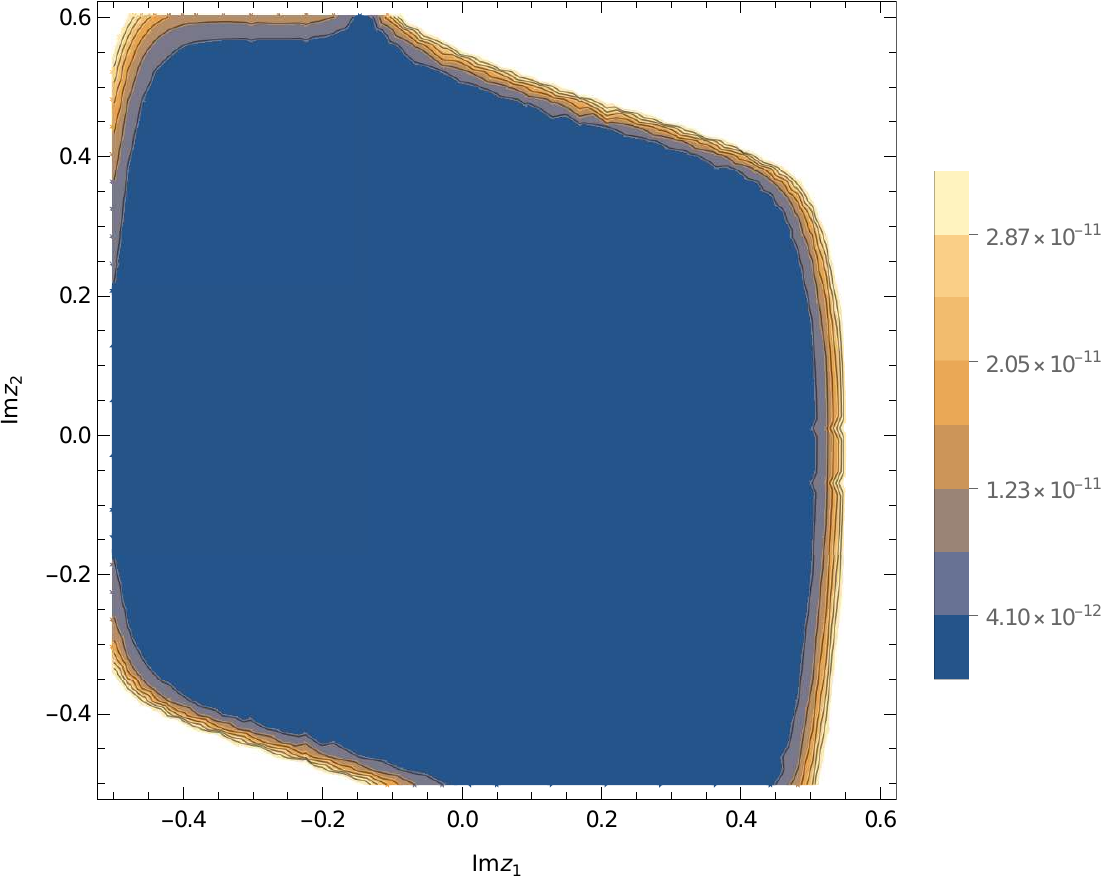}
\caption{}
\end{subfigure}
\caption{The plots where the 10th approximant of (\ref{eq:21})
guarantees certain truncation error bounds for
(\ref{eq:19})}\label{fig:2}
\end{figure}

\begin{table}[h]
\begin{center}
\begin{minipage}{213pt}
\caption{Relative error of {5th} partial sum and {5th} approximant}\label{tab:3}%
\begin{tabular}{@{}llllllll@{}}
\hline
$(z_1,z_2)$&(\ref{eq:19})&(\ref{eq:20})&(\ref{eq:21})\\\hline
$(-0.8,-0.7)$&$-1.1185$&$3.7401\times10^{-01}$&$1.0841\times10^{-04}$\\
$(-0.1,-0.1)$&$-0.1984$&$1.3790\times10^{-09}$&$1.7449\times10^{-12}$\\
$(0.5,-0.7)$&$-0.3396$&$2.0795\times10^{-03}$&$1.4715\times10^{-03}$\\
$(-0.9,0.1)$&$-0.6253$&$2.7472\times10^{-02}$&$1.3940\times10^{-04}$\\
$(0.2,0.3)$&$0.4734$&$2.2374\times10^{-05}$&$1.7797\times10^{-08}$\\
$(0.1,0.8)$&$0.7373$&$2.0631\times10^{-02}$&$3.5972\times10^{-05}$\\
$(0.9,0.9)$&$1.2297$&$2.4591\times10^{+00}$&$2.8455\times10^{-04}$\\
$(2,4)$&$1.7422$&$2.9054\times10^{+05}$&$2.3425\times10^{-02}$\\
$(5,10)$&$1.9697$&$1.0147\times10^{+09}$&$0.3417\times10^{-01}$\\
$(-8,10)$&$-2.0853$&$2.0193\times10^{+05}$&$0.9356\times10^{-01}$\\\hline
\end{tabular}
\end{minipage}
\end{center}
\end{table}

In Figure~\ref{fig:1} (a)--(b) we can see the so-called ``fork
property'' for a branched continued fraction with positive
elements (see, \cite[p.~29]{bod86}). That is, the plots of the
values of even (odd) approximations of (\ref{eq:21}) approach from
below (above) to the plot of the function (\ref{eq:19}).
Figure~\ref{fig:2} (a)--(d) shows the plots, where the 10th
approximant of (\ref{eq:21}) guarantees certain truncation error
bounds for function (\ref{eq:19}).

The numerical illustration of (\ref{eq:20}) and (\ref{eq:21}) is
given in the Table~\ref{tab:3}. Here we can see that the 5th
approximant of (\ref{eq:21}) is eventually a better approximation
to (\ref{eq:19}) than the 5th partial sum of (\ref{eq:20}) is.

Finally, consider the following function of two variables
\begin{align}\notag
\Psi_1(z_1,z_2)&=\psi_1(z_1)+\psi_1(z_2+\psi_1(z_1))=\int_0^\infty\frac{te^{-tz_1}}{1-e^{-t}}dt\\\label{eq:22}
&\quad+\int_0^\infty\frac{s}{1-e^{-s}}\exp\left\{-sz_2-s\int_0^\infty\frac{te^{-tz_1}}{1-e^{-t}}dt\right\}ds,
\end{align}
where $\psi_1(.)$ is trigamma function (see \cite[p.~260]{as64}).

Using the asymptotic expansion for $\psi_1(.)$ given in
\cite[p.~260]{as64}, we find the asymptotic representation for
(\ref{eq:22}) as a formal double Laurent series
\begin{align}\notag
\Psi_1(z_1,z_2)&\approx\sum_{k=0}^\infty\frac{B_k^+}{z_1^{k+1}}+
\sum_{r=0}^\infty\frac{B_r^+}{z_2^{r+1}}
\left(\sum_{s=0}^\infty\left(\sum_{k=0}^\infty\frac{-B_k^+}{z_1^{k+1}z_2}\right)^s\right)^{r+1},\\
&\quad z_i\to\infty,\quad|\arg(z_i)|<\pi,\quad i=1,2,\label{eq:23}
\end{align}
where
\[B_k^+=1-\sum_{r=0}^{k-1}\binom{k}{r}\frac{B_r^+}{k-r+1},\quad k\ge0,\]
are the Bernoulli numbers. Then by Theorem~\ref{th:2} using the
algorithm from Section~\ref{sec:3} we obtain the corresponding
two-dimensional \emph{J}-fraction with independent variables
\begin{equation}\label{eq:24}
\sum_{i_1=1}^2\frac{p_{e_{i(1)}}}{q_{e_{i(1)}}+z_{i_1}}{\atop+}
\sum_{i_2=1}^{i_1}\frac{p_{e_{i(2)}}}{q_{e_{i(2)}}+z_{i_2}}{\atop+}
\sum_{i_3=1}^{i_2}\frac{p_{e_{i(3)}}}{q_{e_{i(3)}}+z_{i_3}}{\atop+}\ldots,
\end{equation}
where
\begin{align*}
p_{e_1+re_2}&=p_{e_2}=1,\quad r\ge0,\\
p_{ke_1+re_2}&=p_{ke_2}=\frac{(k-1)^4}{4(2k-3)(2k-1)},\quad k\ge2,\;r\ge0,\\
q_{ke_1+re_2}&=q_{ke_2}=-1/2,\quad k\ge1,\;r\ge0.
\end{align*}

\begin{figure}[h]%
\centering
\begin{subfigure}[h]{.41\linewidth}
\centering
\includegraphics[width=\linewidth]{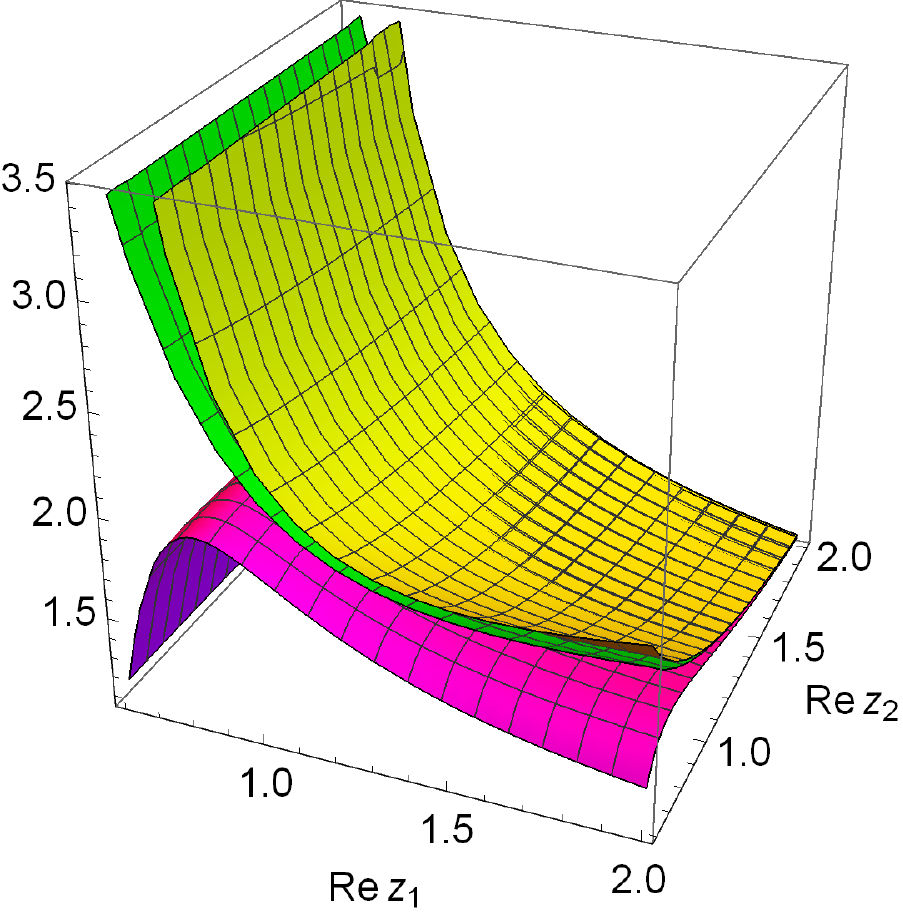}
\caption{\textcolor{rdl}{$\blacksquare$} -- 2nd,
 \textcolor{grl}{$\blacksquare$}  -- $\eqref{eq:22},$
\textcolor{brl}{$\blacksquare$} -- 3rd}
\end{subfigure}\quad\quad
\begin{subfigure}[h]{.37\linewidth}
\centering
\includegraphics[width=\linewidth]{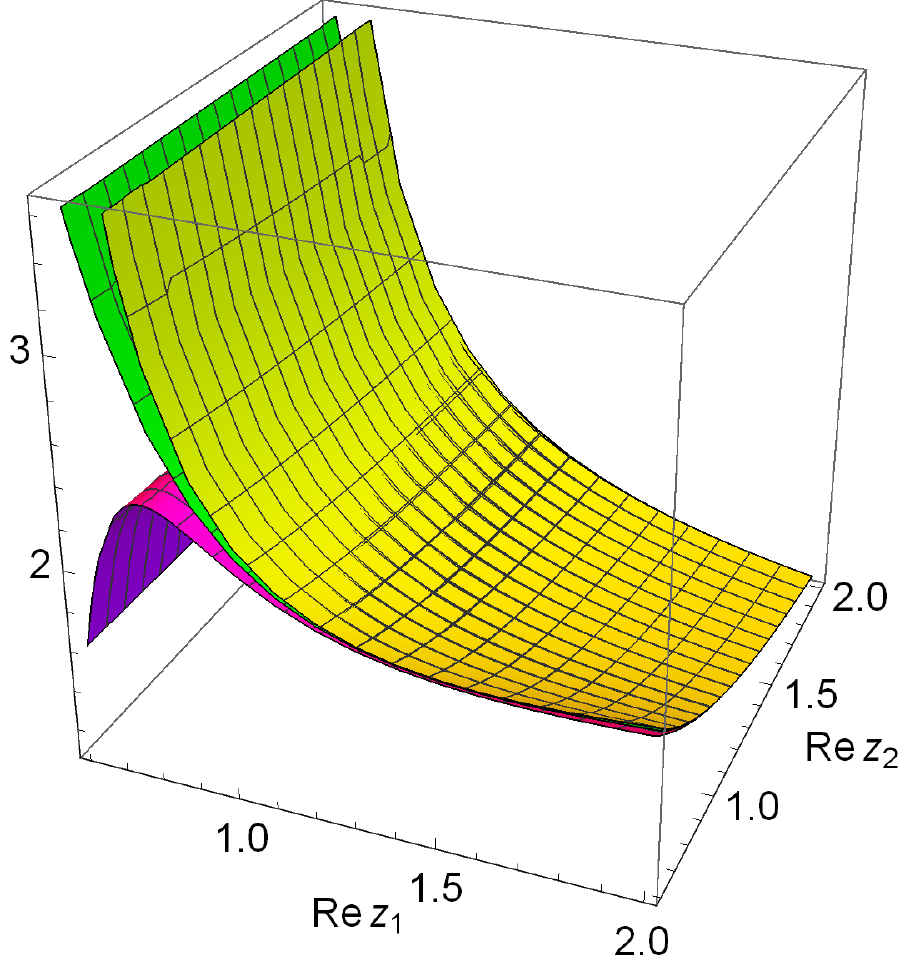}
\caption{\textcolor{rdl}{$\blacksquare$} -- 4th,
 \textcolor{grl}{$\blacksquare$}  -- $\eqref{eq:22},$
\textcolor{brl}{$\blacksquare$} -- 5th}
\end{subfigure}
\caption{The plots of values of the $n$th approximants of
\eqref{eq:24}}\label{fig:3}
\end{figure}

\begin{figure}[h]%
\centering
\begin{subfigure}[h]{.45\linewidth}
\centering
\includegraphics[width=\linewidth]{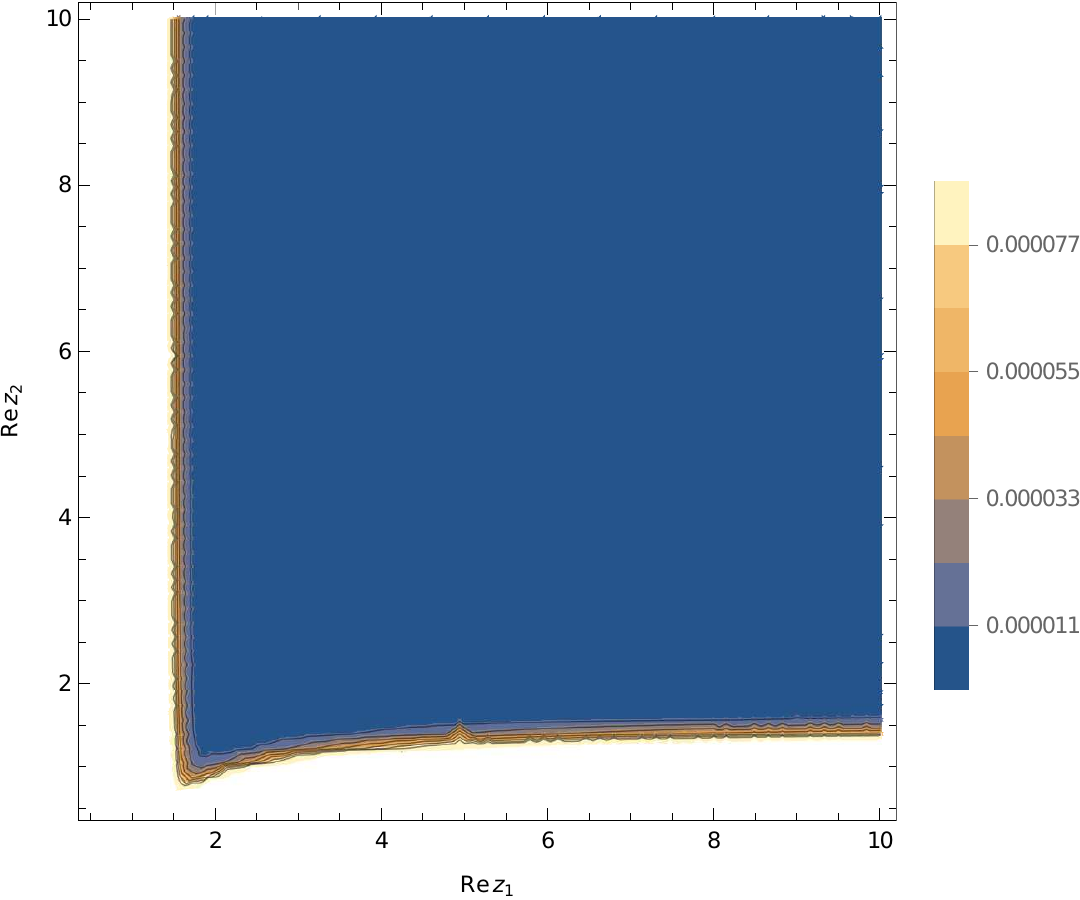}
\caption{}
\end{subfigure}\quad
\begin{subfigure}[h]{.45\linewidth}
\centering
\includegraphics[width=\linewidth]{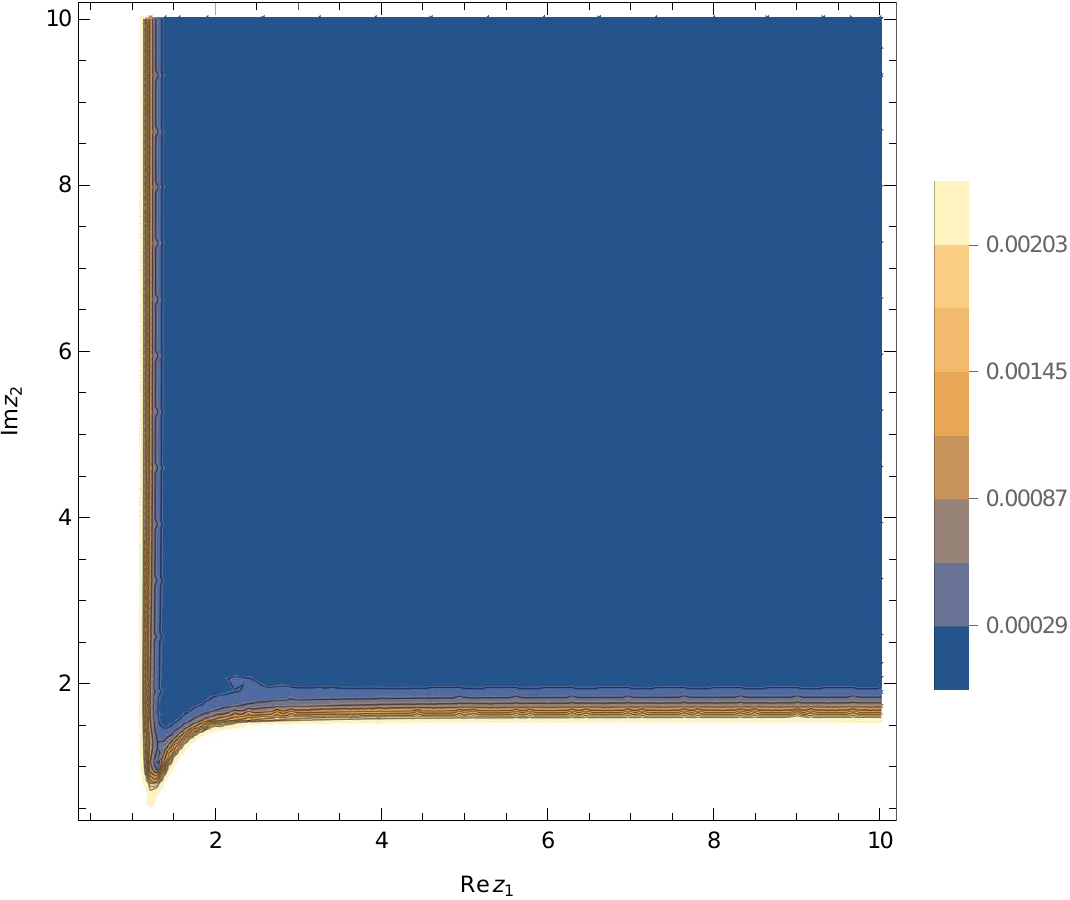}
\caption{}
\end{subfigure}\\[2ex]
\begin{subfigure}[h]{.45\linewidth}
\centering
\includegraphics[width=\linewidth]{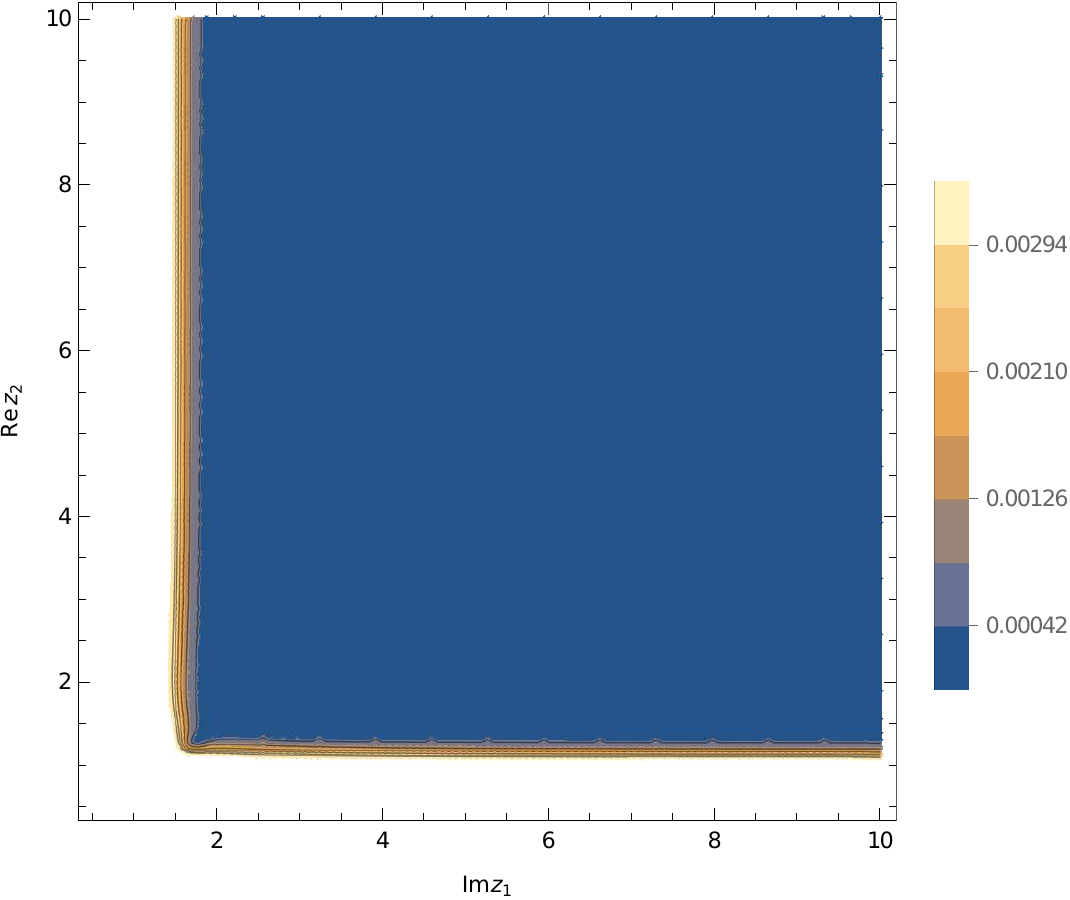}
\caption{}
\end{subfigure}\quad\quad
\begin{subfigure}[h]{.45\linewidth}
\centering
\includegraphics[width=\linewidth]{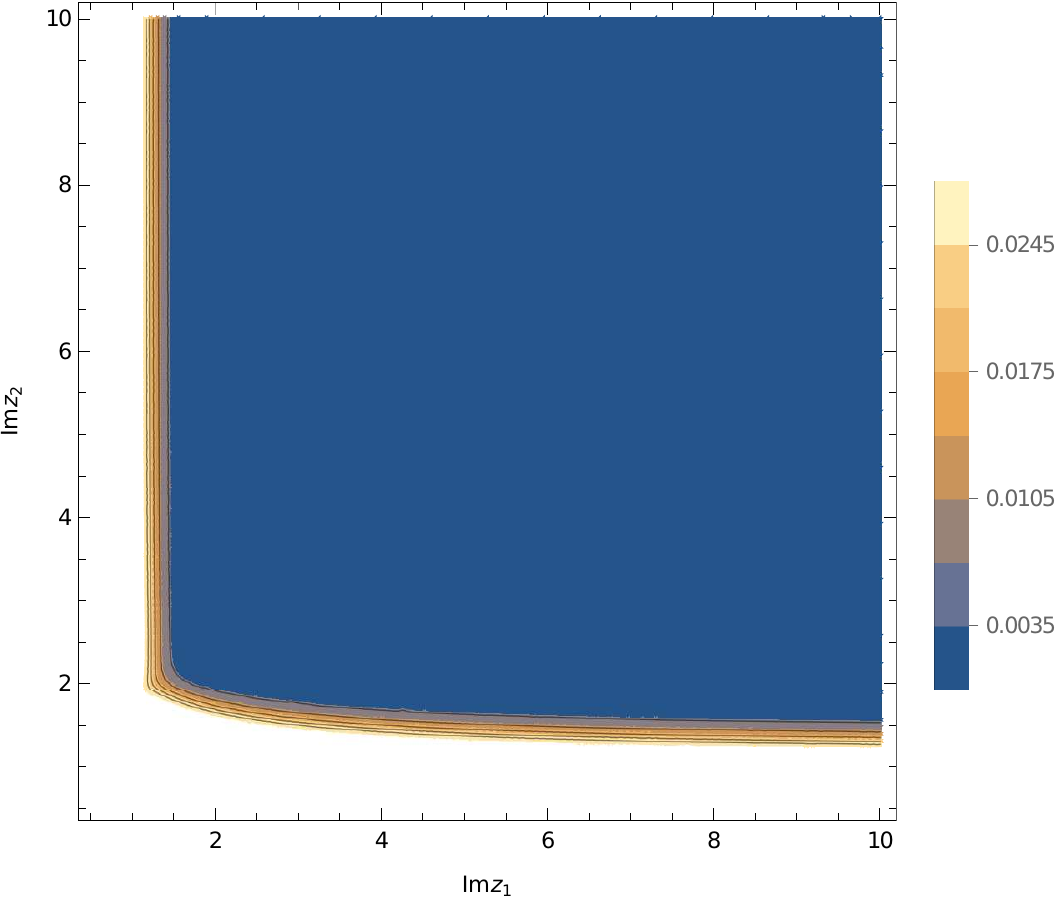}
\caption{}
\end{subfigure}
\caption{The plots where the 10th approximant of (\ref{eq:24})
guarantees certain truncation error bounds for
(\ref{eq:22})}\label{fig:4}
\end{figure}

\begin{table}[h]
\begin{center}
\begin{minipage}{229pt}
\caption{Relative error of {5th} partial sum and {5th} approximant}\label{tab:4}%
\begin{tabular}{@{}llllllll@{}}
\hline
$(z_1,z_2)$&(\ref{eq:22})&(\ref{eq:23})&(\ref{eq:24})\\\hline
$(0.6,0.6)$&$3.9023\times10^{+00}$&$1.1536\times10^{+02}$&$9.1384\times10^{-01}$\\
$(0.9,0.8)$&$2.3654\times10^{+00}$&$8.0419\times10^{+00}$&$4.5692\times10^{-02}$\\
$(1.5,1.4)$&$1.4675\times10^{+00}$&$5.5267\times10^{-02}$&$7.2177\times10^{-04}$\\
$(2,3)$&$9.6031\times10^{-01}$&$1.1374\times10^{-04}$&$4.1428\times10^{-05}$\\
$(10.9)$&$2.2125\times10^{-01}$&$5.8012\times10^{-10}$&$1.0153\times10^{-11}$\\
$(20,40)$&$7.6553\times10^{-01}$&$7.1868\times10^{-14}$&$1.9746\times10^{-15}$\\
$(50,70)$&$3.4585\times10^{-02}$&$8.8000\times10^{-17}$&$1.8847\times10^{-19}$\\
$(100,110)$&$2.0099\times10^{-02}$&$7.8850\times10^{-19}$&$3.6751\times10^{-22}$\\
$(500,1000)$&$3.0025\times10^{-03}$&$2.1889\times10^{-26}$&$1.6536\times10^{-29}$\\\hline
\end{tabular}
\end{minipage}
\end{center}
\end{table}

Plots of the values of the $n$th approximants of the
two-dimensional \emph{J}-fraction with independent variables
(\ref{eq:24}) for function (\ref{eq:22}) are shown in
Figures~\ref{fig:3} (a)--(b). Figure~\ref{fig:4}~(a)--(d) shows
the plots, where the $10$th approximant of (\ref{eq:24})
guarantees certain truncation error bounds for (\ref{eq:22}). The
numerical illustration of (\ref{eq:23}) and (\ref{eq:24}) is given
in the Table~\ref{tab:4}. Here we have results like to the results
in the previous example.




\begin{thebibliography}{999}

\bibitem{as64} Abramowitz, M., Stegun, I.A.: Handbook of Mathematical Functions
With Formulas, Graphs and Mathematical Tables. U.S. Government
Printing Office, NBS, Washington, D.C. (1964)


\bibitem{antdm2020a} Antonova, T.M., Dmytryshyn, R.I.: Truncation error bounds
for branched continued fraction
$\sum_{i_1=1}^N\frac{a_{i(1)}}{1}{\atop+}\sum_{i_2=1}^{i_1}\frac{a_{i(2)}}{1}{\atop+}\sum_{i_3=1}^{i_2}\frac{a_{i(3)}}{1}{\atop+}\ldots.$
Ukr. Math. J. \textbf{72}(7), 1018--1029 (2020)

\bibitem{antdm2020b} Antonova, T.M., Dmytryshyn, R.I.: Truncation error
bounds for branched continued fraction whose partial denominators
are equal to unity. Mat. Stud. \textbf{54}(1), 3--14 (2020)



\bibitem{bar14} Baran,~O.E.: Approximation of Functions of Several Variables by
Branched Continued Fractions with Independent Variables.
Candidate-Degree Thesis (Physics and Mathematics), Pidstryhach
IPPMM NASU, Lviv (2014) (in Ukrainian)


\bibitem{bb21} Bodnar, D.I., Bilanyk, I.B.: Parabolic convergence regions of
branched continued fractions of the special form. Carpathian Math.
Publ. \textbf{13}(1), 619--630 (2021)

\bibitem{bb20} Bodnar, D.I., Bilanyk, I.B.: On the convergence of branched
continued fractions of a special form in angular domains. J. Math.
Sci. \textbf{246}(2), 188--200 (2020)

\bibitem{bod86} Bodnar,~D.I.: Branched Continued Fractions. Naukova
Dumka, Kyiv (1986) (in Russian)

\bibitem{bodm19} Bodnar,~D.I., Dmytryshyn,~R.I.: Multidimensional
associated fractions with independent variables and multiple power
series. Ukr. Math. J. \textbf{71}(3), 370--386 (2019)

\bibitem{bod03} Bodnar, D.I., Hoyenko, N.P.: Approximation of the ratio of
Lauricella functions by a branched continued fraction. Mat. Stud.
\textbf{20}(2), 210--214 (2003) (in Ukrainian)

\bibitem{bod76} Bodnar,~D.I.: Investigation of the convergence of one class of
branched continued fractions. In: Skorobogatko, V.Ya. (ed.)
Continued Fractions and Their Applications, pp. 41--44. Inst.
Math., Acad. Sci. of USSR, Kyiv (1976) (in Russian)


\bibitem{bod98} Bodnar, D.I.: Multidimensional \emph{C}-fractions.
J. Math. Sci. \textbf{90}(5), 2352--2359 (1998)


\bibitem{cl81} Claessens, G.: A generalization of the QD algorithm. J. Comp. and Appl. Math.
\textbf{7}(4), 237--247 (1981)

\bibitem{cuyt88} Cuyt, A.: A multivariate \emph{qd}-like algorithm. BIT \textbf{28}(1), 98--112 (1988)

\bibitem{cuyt83} Cuyt, A.: The \emph{QD}-algorithm and multivariate Pad\'e approximants.
Numer. Math. \textbf{42}(3), 259--269 (1983)

\bibitem{cuve88} Cuyt, A., Verdonk, B.: A review of branched continued fraction
theory for the construction of multivariate rational approximants.
Appl. Numer. Math.  \textbf{4}, 263--271 (1988)


\bibitem{dm17} Dmytryshyn,~R.I.: Convergence of some branched
continued fractions with independent variables. Mat. Stud.
\textbf{47}(2), 150--159 (2017)

\bibitem{dm21cmfth} Dmytryshyn,~R.I.: Convergence of multidimensional \emph{A}- and \emph{J}-fractions
    with independent variables. Comput. Methods Funct.
    \textbf{22}(2), 229--242 (2022)

\bibitem{dm20prsesa} Dmytryshyn,~R.I.: Multidimensional regular
\emph{C}-fraction with independent variables corresponding to
formal multiple power series. Proc. Roy. Soc. Edinburgh Sect. A
\textbf{150}(4), 1853--1870 (2020)

\bibitem{dm05} Dmytryshyn, R.I.: The multidimensional generalization of
\emph{g}-fractions and their application. J. Comp. and Appl. Math.
\textbf{164}--\textbf{165}, 265--284 (2004)

\bibitem{dm12jms} Dmytryshyn,~R.I.: On the expansion of some
functions in a two-dimensional \emph{g}-fraction with independent
variables. J. Math. Sci. \textbf{181}(3), 320--327 (2012)


\bibitem{dm12jat} Dmytryshyn,~R.I.: The two-dimensional
\emph{g}-fraction with independent variables for double power
series. J. Approx. Theory \textbf{164}(12), 1520--1539 (2012)

\bibitem{dm15jms} Dmytryshyn,~R.I.: Two-dimensional generalization
of the Rutishauser \emph{qd}-algorithm. J. Math. Sci.
\textbf{208}(3), 301--309 (2015)

\bibitem{dmsh21} Dmytryshyn, R.I., Sharyn, S.V.: Approximation of functions of several variables by multidimensional
\emph{S}-fractions with independent variables. Carpathian Math.
Publ. \textbf{13}(3), 592--607 (2021)


\bibitem{gr74} Gragg, W. B.: Matrix interpretations and applications of the
continued fraction algorithm. Rocky Mountain J. Math.
\textbf{4}(2), 213--225 (1974)

\bibitem{hers20} Herschel,~J.~F.~W.: A Collection of
Examples of the Applications of the Calculus of Finite
Differences. Printed by J. Smith and sold by J.~Deighton \& Sons,
Cambridge (1820).


\bibitem{jt80} Jones,~W.B., Thron,~W.J.: Continued Fractions:
Analytic Theory and Applications. Addison-Wesley Pub. Co.,
Reading, Mass. (1980)



\bibitem{mudo77} Murphy, J.A., O'Donohoe, M.R.: A class of algorithms for
obtaining rational approximants to functions defined by power
series. J. Appl. Math. Phys. \textbf{28}, 1121-–1131 (1977) 




\bibitem{rut54a} Rutishauser,~H.: Anwendungen des
quotienten-differenzen-algorithmus. Z. Angew. Math. Phys.
\textbf{5}(6), 496--508. (1954)

\bibitem{rut54b} Rutishauser,~H.: Der
quotienten-differenzen-algorithmus. Z. Angew. Math. Phys.
\textbf{5}(3), 233--251 (1954)

\bibitem{rut54c} Rutishauser,~H.: Ein infinitesimales analogon zum
quo\-tien\-ten-diffe\-ren\-zen-algorithmus. Arch. Math.
\textbf{5}(1--3), 132--137 (1954)

\bibitem{siem80} Siemaszko, W.: Branched continued fractions for double power
series. J. Comput. Appl. Math. \textbf{6}(2), 121--125 (1980)



\end{thebibliography}
\end{document}